\documentclass[12pt,dvips]{amsart}
\usepackage{euler, amsfonts, amssymb, latexsym, epsfig,epic,color}
\renewcommand{\familydefault}{ppl}

\newcommand\iso{{\,\cong\,}}
\newcommand\tensor{{\otimes}}

\newtheorem{Theorem}{Theorem} 
\newtheorem{Proposition}{Proposition} 
\newtheorem{Lemma}{Lemma}

\newtheorem{Corollary}{Corollary}
\newtheorem*{Proposition*}{Proposition} 
\newtheorem*{Corollary*}{Corollary}
\newtheorem*{Definition}{Definition} 
\newtheorem*{Theorem*}{Theorem}
\newtheorem*{Lemma*}{Lemma}

\newcommand\onto{\mathop{\to\!\!\!\!\!\to}}
\newcommand\into{\operatorname*{\hookrightarrow}}
\newcommand\Ex{{\em Example. }}

\newcommand\inv{{\rm inv}}

\newcommand\CP{{\mathbb C \mathbb P}}

\newcommand\complexes{{\mathbb C}}
\newcommand\integers{{\mathbb Z}}
\newcommand\rationals{{\mathbb Q}}

\newcommand\Grkn{{{\rm Gr}_k(\complexes^n)}}
\newcommand\Gr{{\rm Gr}}
\newcommand\Oplus{\bigoplus}

\theoremstyle{plain}

\newcommand\dfn{\bf} 

\newcommand\GLn{{{GL_n(\complexes)}}}

\begin{document}
\pagestyle{plain}

\title{Puzzles and (equivariant) cohomology of Grassmannians}
\author{Allen Knutson}
\email{allenk@math.berkeley.edu}
\thanks{AK is supported by the NSF and the Sloan foundation.}
\address{Mathematics Department\\ UC Berkeley\\ Berkeley, California}
\author{Terence Tao}
\thanks{TT is supported by the Clay Mathematics Institute and by
a grant from the Packard foundation.}
\email{tao@math.ucla.edu}
\address{Mathematics Department\\ UCLA\\ Los Angeles, California}
\date{\today}

\begin{abstract}
  The product of two Schubert cohomology classes on a Grassmannian $\Grkn$
  has long been known to be a positive combination of other Schubert classes,
  and many manifestly positive formulae are now available
  for computing such a product (e.g. the Littlewood-Richardson rule,
  or the more symmetric puzzle rule from \cite{Hon2}). 
  Recently in \cite{G} it was shown, nonconstructively,
  that a similar positivity statement holds for {\em $T$-equivariant}
  cohomology (where the coefficients are polynomials). 
  We give the first manifestly positive formula for these
  coefficients, in terms of puzzles using an ``equivariant puzzle piece''. 
  
  The proof of the formula is mostly combinatorial, but requires no
  prior combinatorics, and only a modicum of equivariant cohomology
  (which we include).  As a by-product the argument gives a new proof of
  the puzzle (or Littlewood-Richardson) rule in the
  ordinary-cohomology case, but this proof requires the equivariant
  generalization in an essential way, as it inducts backwards from the
  ``most equivariant'' case.

  This formula is closely related to the one in \cite{MS} for multiplying
  factorial Schur functions in three sets of variables, 
  although their rule does not give a positive formula in the sense of
  \cite{G}. We include a cohomological interpretation of this 
  problem, and a puzzle formulation for it.
\end{abstract}

\maketitle
\tableofcontents

\section{Summary of results}\label{sec:intro}
\newcommand\dv{{\rm div}}
\newcommand\id{{\rm id}}
\newcommand\covers{\to}
\newcommand\clmn{c_{\lambda \mu}^\nu}
\newcommand\PP{{\mathcal P}}

In \cite{Hon2} we introduced a new rule (the puzzle rule) for
computing Schubert calculus (intersection theory on Grassmannians
$\Grkn$), and proved it by reduction to the honeycombs of \cite{Hon1}.
This reduction implicitly involved the somewhat tricky relation
between this cohomology ring and the representation ring of the
general linear group $\GLn$, and so our derivation of the puzzle rule
was somewhat indirect.  In this paper we give an independent and
nearly self-contained proof of the puzzle rule.  The proof is mainly
combinatorial; the only non-combinatorial aspects are a small amount
of equivariant cohomology theory and the Pieri rule $ S_\dv S_\lambda
= \sum_{\lambda': \lambda' \covers \lambda} S_{\lambda'}.$ (For
completeness, we include a combinatorial proof of the Pieri rule in an appendix.)
In particular, we avoid any use of the Littlewood-Richardson
rule\footnote{Indeed, one can use the results in this paper, together
  with the correspondences in \cite{Hon2} and \cite{Bu2}, to prove
  that the Littlewood-Richardson rule computes Schubert calculus,
  though this is not the most direct derivation of this fact.}.  

In the course of our argument we also give a
formula for equivariant Schubert calculus on Grassmannians,
that is manifestly positive in the sense of \cite{G} 
(which shall be defined in a moment);
to our knowledge this is the first such formula.

\subsection{The puzzle rule for $H^*(\Grkn)$.}

We first recall the definition of Schubert calculus, and the puzzle
rule from \cite{Hon2} which computes this calculus.

Let $0 \leq k \leq n$ be fixed integers.  Abusing notation, we let
$n\choose k$ denote the {\em set} of strings $\lambda = \lambda_1
\ldots \lambda_n$ consisting of $k$ ones and $n-k$ zeroes in arbitrary
order, e.g. ${3 \choose 2} = \{ 110, 101, 011 \}$.  This set $n
\choose k$ has an obvious left action\footnote{Readers familiar with
  Schubert classes on flag manifolds may wish to think of $n \choose
  k$ as $S_n$ quotiented by the right-action of $S_{n-k} \times S_k$.
  There are many suggestive hints that the arguments in this paper
  should extend from Grassmannians to flag manifolds, but we have so
  far been unable to extend them.} of the permutation group $S_n$.  In
particular, the long word $w_0$ acts on $n \choose k$ by reversal,
e.g. $w_0\cdot 01101 = 10110$.  If $\lambda \in {n \choose k}$, we define
an {\dfn inversion} of $\lambda$ to be a pair $1 \leq i<j \leq n$ with
$1 = \lambda_i > \lambda_j = 0$.  We denote the set of inversions by
$\inv(\lambda)$ and the number of inversions by $l(\lambda) =
|\inv(\lambda)|$.  Observe that in $n \choose k$ there is a unique
string $\id := 0^{n-k} 1^k$ with no inversions, a unique string $\dv
:= 0^{n-k-1} 10 1^{k-1}$ (assuming $0<k<n$) with one inversion, and a
unique string $w_0\cdot \id := 1^k 0^{n-k}$ with the maximal number $k(n-k)
= \dim_\complexes(\Grkn)$ of inversions.

If $\lambda \in {n\choose k}$ is a string, we let $\complexes^\lambda
:= \Oplus_{i=1}^n \complexes^{\lambda_i}$ denote the corresponding
coordinate $k$-plane in $\complexes^n$, and let $X_\lambda$ be the
{\dfn Schubert cycle} in $\Grkn$ defined as
$$ X_\lambda := \big\{ V_k \in \Grkn \quad: \quad
        \dim (V_k \cap F_i) \geq \dim (\complexes^\lambda \cap F_i),\quad
        \forall i\in [1,n] \big\} $$
where $F_i := \complexes^{0^{n-i} 1^i}$ is the {\em anti}-standard $i$-plane.
Equivalently, $X_\lambda$ is the closure of the set 
$$ \{ V_k \in \Grkn: \lambda_{i} 
        = \dim ((V_k \cap F_i) / (V_k \cap F_{i-1}))
         \hbox{ for all } i = 1, \ldots, n\}.$$
The {\dfn Schubert class} $S_\lambda \in H^*(\Grkn)$ 
is the Poincar\'e dual of the cycle $X_\lambda$.
In particular the degree of $S_\lambda$ is $2l(\lambda)$. These
classes are well-known to give a basis (over $\integers$) for the
cohomology ring $H^*(\Grkn)$, and as such we can expand uniquely the
product $S_\lambda S_\mu$ of any two classes as a sum over the basis
$\{S_\nu\}$, weighted by the structure constants $c_{\lambda\mu}^\nu$
of the multiplication.  These integers $c_{\lambda\mu}^\nu$ are the
concern of (ordinary) ``Schubert calculus''.

Schubert calculus can be computed by many combinatorial rules, most
famously the Littlewood-Richardson rule; we however shall use the more
symmetric puzzle rule from \cite{Hon2}, which we now recall.

Define an (ordinary) {\dfn puzzle piece} as one of the following three
plane figures with labeled edges:
\begin{enumerate}
\item a unit triangle with all edges labeled 0
\item a unit triangle with all edges labeled 1
\item a unit rhombus (two unit triangles glued together along an edge), 
  the two edges clockwise of acute vertices labeled 0, the other two
  labeled 1.
\end{enumerate}
\begin{figure}[htbp]
  \begin{center}
    \leavevmode
    \epsfig{file=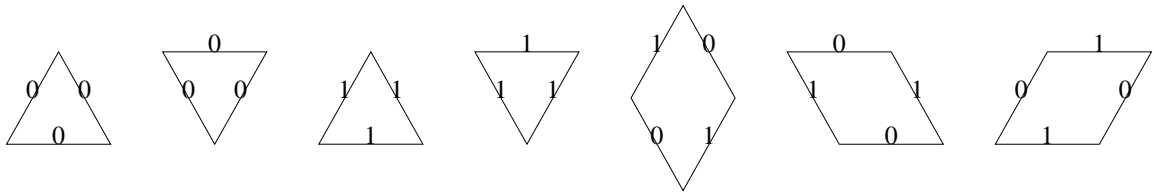,width=6in}
    \caption{The three puzzle pieces, in all their lattice orientations.  
      From left-to-right we have an upward $0$-triangle, a downward 
      $0$-triangle, an upward $1$-triangle, a downward $1$-triangle, 
      a N-S rhombus, a NW-SE rhombus, and a SW-NE rhombus.}
    \label{fig:puzpieces}
  \end{center}
\end{figure}
Note that the set of puzzle pieces is closed under rotation but not reflection
(the reflection of a rhombus puzzle piece is not again a puzzle
piece). See figure \ref{fig:puzpieces}.

\begin{figure}[htbp]
  \begin{center}
    \leavevmode
    \epsfig{file=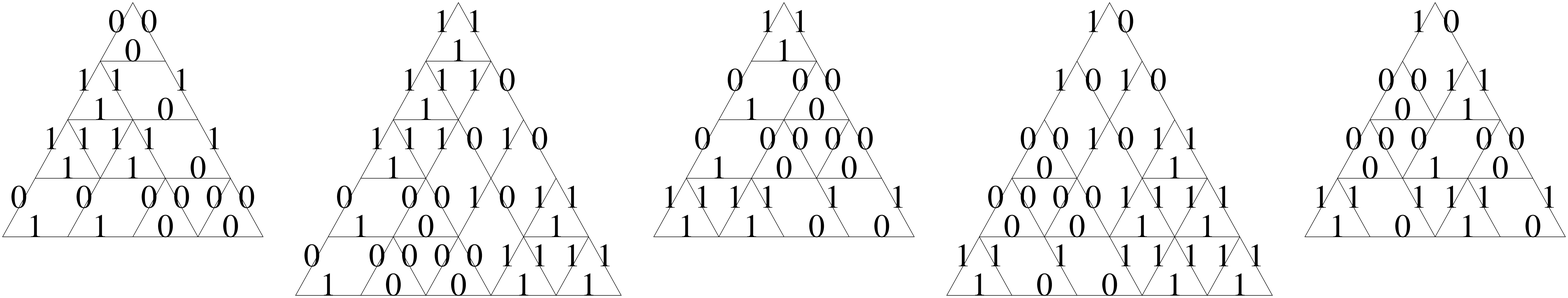,width=6in}
    \caption{Some examples of puzzles.  If $P$ denotes the right-most puzzle, 
      then $\partial P = \Delta_{1001,0101,0101} = \Delta_{1001,0101}^{1010}$.}
    \label{fig:puzex}
  \end{center}
\end{figure}

Define an (ordinary) {\dfn puzzle} as a decomposition of
an equilateral triangle
into triangles and rhombi, all edges labeled 0 or 1, such that each
region is a puzzle piece. (Alternately, one can speak of attaching
puzzle pieces together, with edges required to match up as in a jigsaw
puzzle.) We will always align our puzzles to have a South side,
Northwest\footnote{Our definition of northwest will be at a 
  $60^\circ$ angle to north, rather than $45^\circ$.  Similarly for
  southwest, etc.} side and Northeast side; this forces the edges of
puzzle pieces to be oriented E-W, NW-SE, or NE-SW, the triangles to be
oriented upward or downward, and the rhombi to be oriented N-S, NW-SE,
or SW-NE. Some examples of puzzles are pictured in figure
\ref{fig:puzex}.

Define a {\dfn labeled equilateral triangle} to be an upward-pointing
equilateral triangle of some integer side-length $n$, with the $3n$
unit edges on the boundary labeled either $0$ or $1$.  Clearly every
puzzle $P$ induces a labeled equilateral triangle $\partial P$, which
we refer to as the {\dfn boundary} of $P$.

Given any three strings $\lambda, \mu, \nu \in {n \choose k}$, we let
$\Delta_{\lambda \mu \nu}$ denote the labeled equilateral triangle
with NW side labeled $\lambda$, NE side labeled $\mu$, and S side
labeled $\nu$ (all read clockwise).  We also let 
$\Delta_{\lambda \mu}^\nu = \Delta_{\lambda,\mu,w_0\cdot \nu}$ denote the labeled
equilateral triangle with NW side $\lambda$, NE side $\mu$, and S side
$\nu$, all read left-to-right.  If $P$ is a puzzle, we let $P_{NW}$,
$P_{NE}$, $P_S$ be the three strings of labels of $\partial P$ read
clockwise, thus
$$ \partial P = \Delta_{P_{NW} P_{NE} P_S} 
        = \Delta_{P_{NW} P_{NE}}^{w_0\cdot P_S}.$$
We will also call a puzzle with boundary $\Delta_{\lambda \mu}^\nu$
a {\dfn $\Delta_{\lambda \mu}^\nu$-puzzle}.

Our first main result shall be a new, and essentially self-contained, proof of
the following theorem.

\begin{Theorem}[Puzzles compute Schubert calculus]\label{thm:puzcount}
\cite{Hon2}
  Let $0 \leq k \leq n$, and let $\lambda,\mu,\nu$ be three elements
  of $n \choose k$ indexing Schubert classes $S_\lambda,S_\mu,S_\nu$
  in $H^*(\Grkn)$.  Then the following (equivalent) statements hold:
  \begin{enumerate}
  \item The intersection number $\int_{\Grkn} S_\lambda S_\mu S_\nu$
    is equal to the number of puzzles $P$ with $\partial P =
    \Delta_{\lambda \mu \nu}$.
  \item The structure constant $c_{\lambda\mu}^\nu$ is equal to the
    number of puzzles with $\partial P = \Delta_{\lambda \mu}^\nu$.
  \item $$ S_\lambda\, S_\mu 
  = \sum_{\hbox{puzzles $P$: } \atop {\hbox{$P_{NW} = \lambda, P_{NE} = \mu$}}}
        S_{w_0\cdot P_S} $$
  \end{enumerate}
\end{Theorem}

This first formulation, in terms of Schubert intersection numbers,
realizes several symmetries evident in that problem.  Note that the
$120^\circ$ rotation of a puzzle is again a puzzle, corresponding to
the fact\footnote{One would also expect a similar geometric
  interpretation of the commutativity property $S_\lambda S_\mu =
  S_\mu S_\lambda$; the puzzle rule can indeed be shown to be
  commutative but this turns out to be much more non-trivial.} that
$$\int_{\Grkn} S_\lambda S_\mu S_\nu = \int_{\Grkn} S_\mu S_\nu S_\lambda = 
\int_{\Grkn} S_\nu S_\lambda S_\mu.$$ 
We include here the standard proof that these integrals are a priori positive.
They are visibly computing the number of signed intersection points of
three Schubert cycles, perturbed to be transverse. 
It turns out to be possible to achieve this perturbation by replacing the
standard flag $(F_i)$ by two other generic flags, which means the three
transverse cycles are again complex subvarieties.  
Then the intersection points 
all have positive sign. Unfortunately this simple proof, which
generalizes to arbitrary flag manifolds for arbitrary groups $G$,
does not provide a formula (under most people's notions of ``formula'').

From degree considerations we see that the structure constants $\clmn$
vanish unless $l(\lambda) + l(\mu) = l(\nu)$.  We invite the reader to
see how this simple fact can also be deduced from the puzzle rule.

If $\lambda \in {n \choose k}$, define the {\dfn dual string}
$\lambda^* \in {n \choose n-k}$ to be the string $w_0\cdot \lambda$ with
all $0$s and $1$s exchanged; thus for instance $01100^* = 11001$. The
dual string $\lambda^*$ gives a Schubert class $S_{\lambda^*}$ on the
dual Grassmannian $\Gr_{n-k}(\complexes^n)$.  
Similarly, given a puzzle $P$, we can define the {\dfn dual puzzle}
$P^*$ by reflecting $P$ left-right and exchanging $1$s and $0$s
everywhere.  Observe that this gives a puzzle-theoretic proof of the equation
$$ \int_{\Grkn} S_\lambda S_\mu S_\nu 
        = \int_{\Gr_{n-k}(\complexes^n)} S_{\mu^*} S_{\lambda^*} S_{\nu^*}, $$

Grassmann duality $\Grkn \,\iso\, \Gr_{n-k}(\complexes^n)$
gives a geometric proof of this identity, as follows.
It takes a $k$-plane to its perpendicular $(n-k)$-plane (with respect to
the standard Hermitian form on $\complexes^n$), 
and the Schubert variety $X_\lambda$
to the {\em opposite} Schubert variety $w_0\cdot X_{\lambda^*}$
(thinking of $w_0$ as the antidiagonal permutation matrix).
Since the transformation $w_0$ is deformable to the identity transformation,
$w_0\cdot X_{\lambda^*}$ again represents the Schubert class $S_{\lambda^*}$.

\begin{figure}[htbp]
  \begin{center}
    \leavevmode
    \epsfig{file=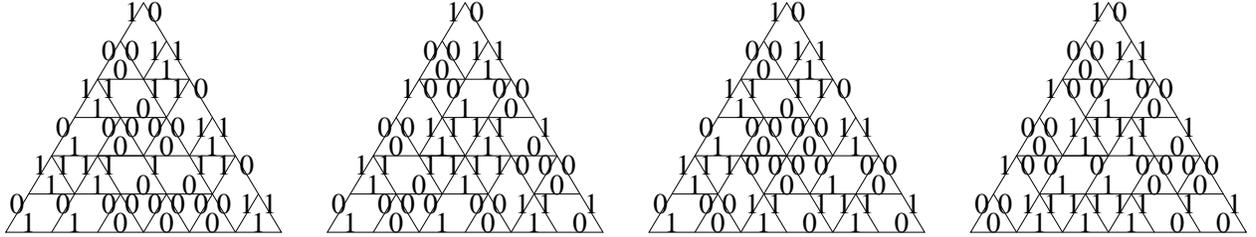,width=6.5in}
    \caption{The four puzzles $P$ with $P_{NW} = P_{NE} = 010101$,
      computing $S_{010101}^2 = S_{110001} + 2\, S_{101010} + S_{011100}$
      in $H^*(\Gr_3(\complexes^6))$.}
    \label{fig:prodex}
  \end{center}
\end{figure}

The third formulation in theorem \ref{thm:puzcount}
is very suitable for computations; an example is in figure
\ref{fig:prodex}.

\newcommand\tS{{\tilde S}}
\subsection{A new puzzle piece, for $H^*_T(\Grkn)$.}

To prove theorem \ref{thm:puzcount} we shall 
generalize the
result so that it computes the answer to a harder question, namely the
product structure in the {\em $T$-equivariant} cohomology of
Grassmannians.
In section \ref{sec:eqvtH} we recall the (very few) necessary facts
about $T$-equivariant cohomology we need to set up this question.
For now, we need only four:
\begin{itemize}
\item the equivariant cohomology ring $H^*_T(\Grkn)$ is naturally
  a graded module over the polynomial ring $\integers[y_1,\ldots,y_n]$
  (itself the equivariant cohomology of a point);
\item $H^*_T(\Grkn)$ has a natural basis of ``equivariant Schubert classes'' 
  $\{ \tS_\lambda \}$, with $\tS_\lambda$ having degree $2l(\lambda)$;
\item there is a natural forgetful map
  $H^*_T(\Grkn) \to H^*(\Grkn)$ to ordinary cohomology, which 
  consists of setting all the $y_i$ to $0$;
\item this forgetful map takes each equivariant Schubert class $\tS_\lambda$ to the corresponding   ordinary Schubert class $S_\lambda$.
\end{itemize}

In particular, one can speak of ``equivariant Schubert calculus'',
which concerns the structure constants $c_{\lambda\mu}^\nu \in
\integers[y_1,\ldots,y_n]$ in the product expansion%
\footnote{In this paper, all summations over Greek indices shall
  range over $n \choose k$.}  
$\tS_\lambda \tS_\mu = \sum_\nu c_{\lambda\mu}^\nu \tS_\nu$.  By
degree considerations in this graded ring, we know $\deg
c_{\lambda\mu}^\nu$ is a homogeneous polynomial of degree $l(\lambda)
+ l(\mu) - l(\nu)$.  In particular, $c_{\lambda\mu}^\nu$ vanishes when
$l(\lambda) + l(\mu) < l(\nu)$, and agrees with the ordinary structure
constants $c_{\lambda\mu}^\nu$ when $l(\lambda) + l(\mu) = l(\nu)$
(which is why we can safely use the same notation for both).

It is not hard to show that the equivariant structure constants 
$c_{\lambda\mu}^\nu$ 
actually live in the subring $\integers[y_2-y_1,\,y_3-y_2,\ldots,y_n-y_{n-1}]$.
In \cite{G} it is proven that written as polynomials in these differences, 
the structure constants have positive integer coefficients
(this was first conjectured by Dale Peterson). 
As in the non-equivariant case, the proof does not directly give a
formula for the $c_{\lambda\mu}^\nu$.

\begin{figure}[htbp]
  \begin{center}
    \epsfig{file=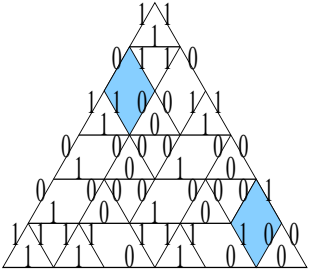, height=1.4in, width=1.6in}
    \caption{A puzzle with two equivariant pieces, which are shaded. 
      The left equivariant piece has weight $y_4-y_1$, the right 
      $y_5-y_4$, so this puzzle contributes $(y_4-y_1)(y_5-y_4)$
      to the calculation of $c_{100101,101010}^{\ \ \ 110100}$.}
    \label{fig:sampleeqvt}
  \end{center}
\end{figure}

To compute these $c_{\lambda\mu}^\nu$, we need to generalize our
notion of puzzle a bit. We introduce the {\dfn equivariant puzzle piece}: 
this is the same as the N-S rhombus puzzle piece but with
the $1$s and $0$s interchanged.  A puzzle using some equivariant pieces%
\footnote{Since equivariant Schubert calculus generalizes ordinary,
  we can safely call these again ``puzzles'' and need not introduce a
  term ``equivariant puzzles''.  Rather, one might call a
  puzzle {\dfn ordinary} if one wanted to emphasize that 
  it happens to contain no equivariant pieces.}  
is given in figure \ref{fig:sampleeqvt}.

\begin{figure}[htbp]
  \begin{center}
    \input{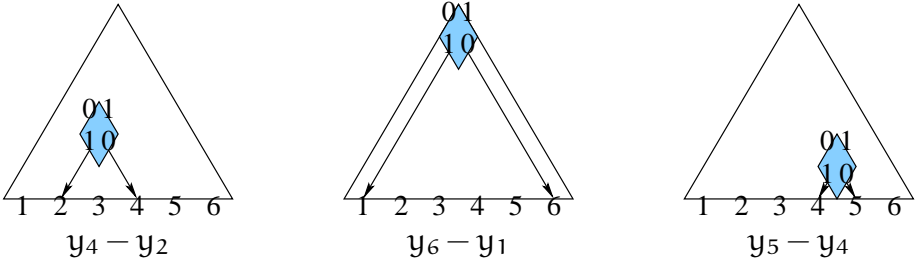}
    \caption{Locations of some equivariant puzzle pieces, 
      and their corresponding weights.}
    \label{fig:piecewts}
  \end{center}
\end{figure}

To each equivariant piece $p$ in a puzzle, we associate a {\dfn weight} 
$wt(p)$, which we compute by dropping lines SW and SE from
the piece until they poke out the $i$th and $j$th place on the South
side and then setting $wt(p) := y_j - y_i$. See figure
\ref{fig:piecewts}. The weights of the pieces in figure
\ref{fig:sampleeqvt} are given as example.  We can then associate a
{\dfn weight} $wt(P)$ to every puzzle $P$ by defining $wt(P) = \prod_p
wt(p)$, where $p$ ranges over the equivariant pieces of $P$. (An empty
product is taken to be $1$, of course.)

Observe that we necessarily have $i < j$ in the above definition of
$wt(p)$.  In particular, the weight $wt(P)$ of a puzzle can be
expressed as a positive combination of $y_2-y_1,\,y_3-y_2,\ldots,y_n-y_{n-1}$.

The main result of this paper is

\begin{Theorem}[
  Puzzles compute equivariant Schubert calculus]\label{thm:eqvtpuzcount}
  Let $0 \leq k \leq n$, and let $\lambda,\mu,\nu$ be three elements
  of $n \choose k$, indexing equivariant Schubert classes
  ${\tS}_\lambda,{\tS}_\mu,{\tS}_\nu$ in $H^*_T(\Grkn)$.  Then the
  following (equivalent) statements hold:
  \begin{enumerate}
  \item The structure constant $\clmn$ is equal to the sum of the
    weights of all puzzles $P$ with $\partial P = \Delta_{\lambda \mu}^\nu$.  
    In particular, we explicitly demonstrate for $\Grkn$
    the abstract positivity result in \cite{G}.
  \item $$ {\tS}_\lambda\, {\tS}_\mu 
  = \sum_{\hbox{puzzles $P$: } \atop \hbox{$P_{NW} = \lambda, P_N = \mu$}}
        wt(P) \, S_{w_0\cdot P_S} $$
  \end{enumerate}
\end{Theorem}

This obviously implies the second and third formulations of theorem
\ref{thm:puzcount}.  
There is no close analogue of the first 
formulation.\footnote{The closest analogue would be to compute
$\int_\Grkn \tS_\lambda \tS_\mu \tS_\nu = c_{\lambda\mu\nu}^{w_0\cdot id}
:= \sum_\rho c_{\lambda\mu}^\rho c_{\rho\nu}^{w_0\cdot id}$. So a
positive formula for the structure constants gives a positive formula
for the integrals. The converse does not seem to be obviously true.}
In ordinary cohomology the three formulations could be equated via the formula 
$\int_{\Grkn} S_\lambda S_\mu = \delta_{\lambda,w_0\cdot \mu}$, 
but this identity does not hold in equivariant cohomology.
In particular, we should not lament the symmetry lost by including the
non-rotatable equivariant piece, since the problem itself is less symmetric. 
On the other hand, the dual $P^*$ of an equivariant puzzle
is still an equivariant puzzle, giving an equality 
$c_{\lambda\mu}^\nu = \overline{c_{\mu^*\lambda^*}^{\nu^*}}$,
where the bar is defined by $\overline{y_i} := -y_{n+1-i}$.
This again follows from Grassmann duality (the two coefficients are not
equal, because $w_0$ is not deformable to the identity through 
$T$-invariant maps).

\begin{figure}[htbp]
  \begin{center}
    \leavevmode
    \epsfig{file=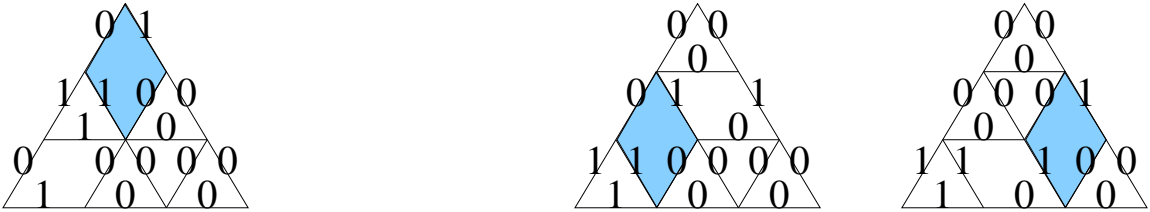,width=6in}
    \caption{The left puzzle computes 
      $\tS_{010} \tS_{100} = (y_3-y_1) \tS_{100}$, the
      right two
      $\tS_{100}\tS_{010} = (y_2-y_1)\tS_{100} +(y_3-y_2)\tS_{100}$.}
    \label{fig:prodexeqvt}
  \end{center}
\end{figure}

We give an example in figure \ref{fig:prodexeqvt}, computing the products
$\tS_{100} \tS_{010}$ and $\tS_{010} \tS_{100}$.
These are of course equal (the ring $H^*_T(\Grkn)$ is commutative),
but this is very nonobvious from the formula.

This paper is organized as follows.  In the ``geometric'' part
(sections \ref{sec:eqvtH} -- \ref{sec:pieri}, plus an optional Appendix) of
the paper we set up the machinery from equivariant cohomology which we
will need, culminating in the equivariant Pieri identities for
Schubert classes and structure constants.  In the ``combinatorial''
part (sections \ref{diag-sec} -- \ref{sec:cscab}) we show that the
equivariant puzzle rule obeys these Pieri identities, which will imply
theorem \ref{thm:eqvtpuzcount} (and hence theorem \ref{thm:puzcount})
by an induction argument.  We remark that this induction only seems to
be available in the equivariant setting, so we cannot give a
completely non-equivariant proof of theorem \ref{thm:puzcount} by
these techniques.  

Finally, in section \ref{sec:ms} we compare the
results here with those in \cite{MS} for multiplying factorial Schur
functions (which are nothing other than double (i.e. equivariant)
Schubert polynomials for Grassmannian permutations). They solve a
different problem, which also reduces to ordinary Schubert calculus
when $l(\nu) = l(\lambda)+l(\mu)$. We introduce cohomological
formulations of their problem, and a reformulation of their rule
in terms of ``MS-puzzles''.

We have had many useful conversations with Chris Woodward, our
coauthor on \cite{Hon2}. Our approach to ordinary Schubert calculus
by inducting from the equivariant counterpart was inspired by \cite{MS}.  
We are very grateful to Anda Degeratu for suggesting the name ``puzzle''.

\section{Equivariant cohomology, especially of Grassmannians}
\label{sec:eqvtH}

Fix $0 \leq k \leq n$.  In this section we give a \emph{combinatorial}
definition of the equivariant cohomology ring $H^*_T(\Grkn)$, which we
interpret as lists of polynomials indexed by $n\choose k$ satisfying
some congruence conditions. We then invoke some standard facts about
equivariant cohomology to determine that this ring is indeed the ring
of equivariant cohomology classes on the Grassmannian, and is equipped
with a basis of equivariant Schubert classes which map to ordinary
Schubert classes under the forgetful map\footnote{We chose this rather
  odd approach to $H^*_T(\Grkn)$ to emphasize the point that one does
  not actually need much equivariant cohomology theory to prove
  Theorem \ref{thm:puzcount}, and one could in fact just think of
  $H^*_T(\Grkn)$ as an abstract ring of lists of polynomials to be
  manipulated combinatorially without ever having to understand what
  the functor $H^*_T$ means.}.  Our reference for combinatorial
properties of equivariant cohomology is \cite{GZ}.

\subsection{A combinatorial description of $H^*_T(\Grkn)$.}

Begin by {\em defining} $H^*_T(pt)$ to be the polynomial ring
$\integers[y_1,\ldots,y_n]$ in $n$ variables (without yet worrying
what ``$H^*_T$'' means in general).  Define $\Oplus_{n\choose k}
H^*_T(pt)$ to be the space of all lists of polynomials $\alpha =
(\alpha|_\lambda)$, indexed by elements $\lambda \in {n\choose k}$.
This is clearly a commutative ring, and a $H^*_T(pt)$-module (where
$H^*_T(pt)$ acts diagonally on each term $\alpha|_\lambda$ of the
list).

Suppose that $\alpha \in \Oplus_{n\choose k} H^*_T(pt)$.
Call $\alpha$ a {\dfn class} if it satisfies the {\em GKM conditions}:
\begin{quote}
  For each pair $\lambda$, $\lambda' \in {n\choose k}$
  differing\footnote{Equivalently, 
    we have $\lambda' = (i \leftrightarrow j) \lambda$, 
    where $(i \leftrightarrow j) \in S_n$ is the transposition of $i$
    and $j$.}  only in places $i$ and $j$, the difference
  $\alpha|_\lambda - \alpha|_{\lambda'}$ should be a multiple of
  $y_i-y_j$.
\end{quote}

\begin{figure}[htbp]
  \begin{center}
    \input{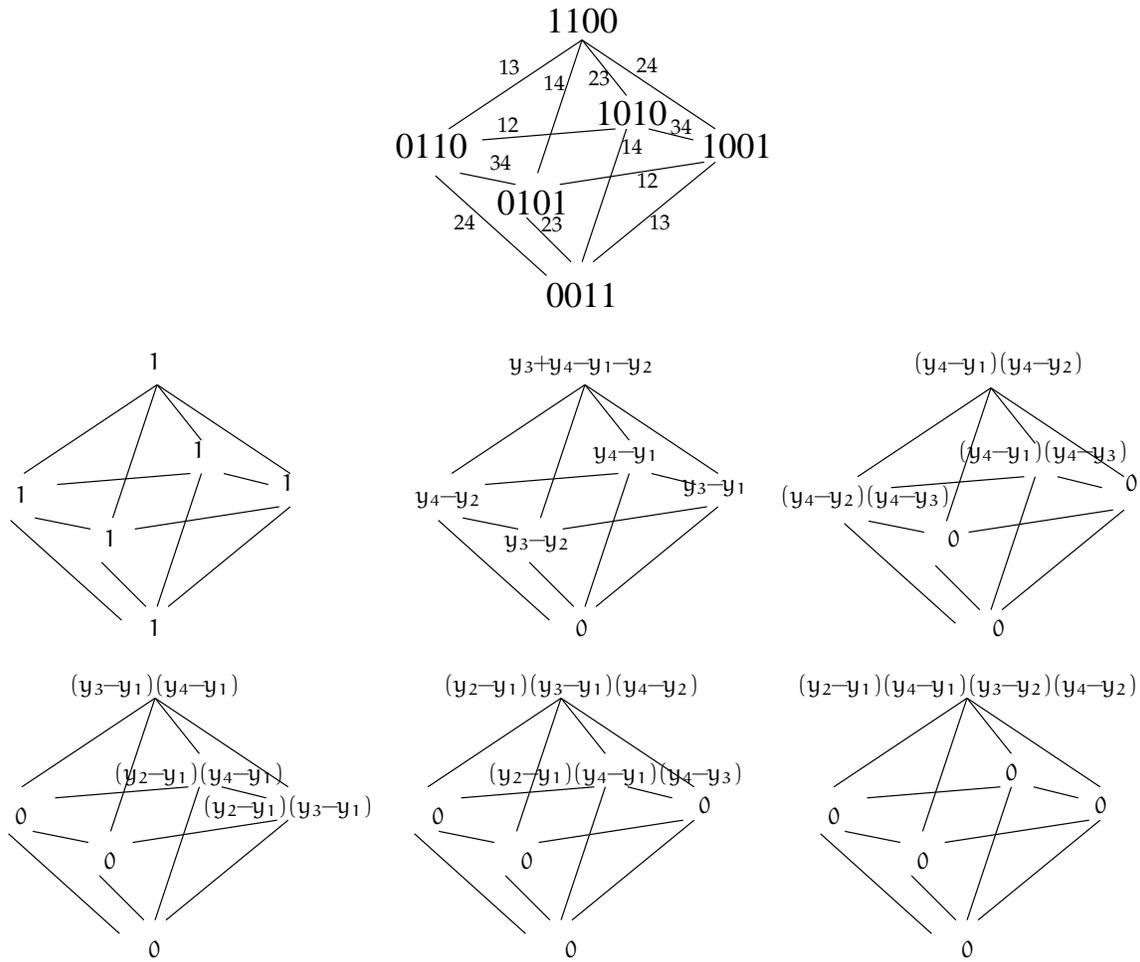}
    \caption{The six Schubert classes in $H^*_T({\rm Gr}_2(\complexes^4))$.
      The top figure is the lattice $4 \choose 2$, with edges labeled
      $ij$ between two strings differing in only the $i,j$ positions,
      thereby causing a GKM condition. }
    \label{fig:Gr24}
  \end{center}
\end{figure}

{\em Examples.} The list $\alpha|_\lambda := 1$ is a class, since
all the relevant differences are $0$.
The list $ \alpha|_\lambda := \sum_{i=1}^n \lambda_i y_i $ is also 
a class\footnote{Topologically, this class $\alpha$ arises as the equivariant
  first Chern class of the $k$th exterior power of the tautological
  $k$-plane bundle on the Grassmannian. 
}, where the multiples are all $1$. 
For each $\mu$, the list
$\alpha|_\lambda := \delta_{\lambda,\mu} \prod_{i<j} (y_i-y_j)$
is also a class. In figure \ref{fig:Gr24} are a list of some very
special classes in $H^*(Gr_2(\complexes^4))$.

Now {\em define} $H^*_T(\Grkn) \subseteq \Oplus_{n\choose k}$ to be
the set of all classes. It is obviously a subring and a submodule of
$\Oplus_{n\choose k} H^*_T(pt)$.

\newcommand\supp{{\rm supp}}
Define the {\dfn support} $\supp(\alpha)$ of a class $\alpha$ as the
set $\lambda \in {n\choose k}$ such that $\alpha|_\lambda \neq 0$.
Recall that $n\choose k$ is a lattice, where the partial order is
given by $\lambda' \geq \lambda$ 
if one has $\sum_{i=1}^j \lambda'_i \geq \sum_{i=1}^j \lambda_i$ for
all $j=1,\ldots,n$.  Say that $\alpha$ is {\dfn supported above $\lambda$} 
if one has $\lambda' \geq \lambda$ for all $\lambda' \in \supp(\alpha)$.

Let $\lambda \in {n \choose k}$, and let $\alpha$ be a class supported
above $\lambda$.  This forces $\alpha|_\lambda$ to be a multiple
of\break $\prod_{(i,j) \in \inv(\lambda)} (y_j-y_i)$.  If we have the
stronger relationship 
$\alpha|_\lambda = \prod_{(i,j) \in \inv(\lambda)} (y_j-y_i)$, 
and also $\alpha|_\mu$ is homogeneous of degree\footnote{In this paper
  we give each generator $y_i$ a degree of 1, although from the
  cohomological considerations below it could be argued that the $y_i$
  really deserve to have degree 2.  We believe however that setting
  $\deg(y_i) = 2$ here would be too confusing.  We remark that if we
  replaced the equivariant cohomology ring with the equivariant Chow
  ring (which is equivalent for $\Grkn$) then the $y_i$ genuinely do
  have degree 1.} 
$l(\lambda)$ for all $\mu \in {n \choose k}$, we call $\alpha$ a 
{\dfn Schubert class corresponding to $\lambda$}. For some examples, see
the Schubert classes in $H^*_T({\rm Gr}_2(\complexes^4))$ in figure
\ref{fig:Gr24}.

\begin{Lemma}[Schubert classes are unique]\label{lem:uniq}  
  For each $\lambda \in {n \choose k}$ there is at most one Schubert
  class $\alpha$ corresponding to $\lambda$.
\end{Lemma}

\begin{proof}
  Suppose for contradiction that there were two distinct Schubert
  classes $\alpha, \alpha'$ corresponding to $\lambda$. Let $\mu$ be a
  minimal element of the support of the class $\alpha - \alpha'$.
  Since $\alpha$, $\alpha'$ agree on $\lambda$ and are supported above
  $\lambda$, we have $\mu > \lambda$.  By the GKM conditions this
  forces $(\alpha-\alpha')|_\mu$ to be a multiple of 
  $\prod_{(i,j) \in \inv(\mu)} (y_j-y_i)$.  But this contradicts the
  fact that $\alpha-\alpha'$ is homogeneous of degree $l(\lambda) < l(\mu)$.
\end{proof}

To prove existence of Schubert classes is a little trickier.  We now
give a topological proof that there exists a Schubert class
$\tS_\lambda$ for each $\lambda \in {n \choose k}$.  We also give
a purely combinatorial proof in the Appendix.

\subsection{$T$-equivariant cohomology.}

\newcommand\TTop{T-{\bf Top}}
Let $T:=(S^1)^n$ be a torus, and $\TTop$ the category of topological
spaces with a continuous $T$-action, the morphisms being equivariant maps.
Then $T$-equivariant cohomology is a contravariant functor 
from $\TTop$ to supercommutative\footnote{Elements of odd degree
anticommute with one another, as in ordinary cohomology. 
In the case of interest, the Grassmannian,
there are none of these anyway and the ring is therefore commutative.} 
graded rings. 
We do not define it here, as all we need are a few of its properties:

\begin{enumerate}
\item $H^*_T(pt)$ is (as promised) 
  a polynomial ring $\integers[y_1,\ldots y_n]$, whose generators are
  given formal degree $2$, and correspond to a basis of the dual of
  the Lie algebra of $T$. In particular we can think of weights of $T$
  as giving elements of $H^2_T(pt)$ 
  (i.e. the weights are linear combinations of the $\{y_i\}$).
\item If $X$ is a $T$-invariant oriented cycle in a compact oriented
  manifold $Y$, then $X$ determines\footnote{ In ordinary cohomology
    $[X]$ can be defined via Poincar\'e duality, but that is not
    available in equivariant cohomology. Nonetheless this should be
    the intuition.} an equivariant cohomology
  class on $Y$, which we will denote $[X]$.  The degree of $[X]$ is
  the codimension of $X$ in $Y$.  If $f \in Y^T$ is not in $X$, then
  the pullback of $[X]$ to $\{f\}$ is zero.
\item There is a natural ``forgetful'' map from $H^*_T(Y) \to H^*(Y)$,
  which factors as
$$H^*_T(Y)\quad  \onto \quad H^*_T(Y) \tensor_{H^*_T(pt)} \integers 
\quad \into \quad H^*(Y).$$
  It takes the equivariant class 
  $[X] \in H^*_T(Y)$ associated to a $T$-invariant cycle to the
  ordinary Poincar\'e dual of $X$ in $H^*(Y)$.
\item Given a $T$-space $Y$, there are two natural equivariant maps 
  associated, $Y^T \into Y \onto pt$, where $Y^T$ is the set of 
  fixed points. These induce ring homomorphisms 
  backwards in equivariant cohomology, 
  $H^*_T(pt) \to H^*_T(Y) \to H^*_T(Y^T) \iso H^*(Y^T) \tensor H^*_T(pt)$. 
  In other words, the functor $H^*_T$ takes values in the
  category of algebras over $H^*_T(pt)$.
\end{enumerate}

Now assume that $Y$ is a smooth projective variety, and $T$ acts on
$Y$ algebraically with isolated fixed points.  (In particular, $Y$
could be the Grassmannian $\Grkn$).  Then much more is true (see \cite{GKM}):

\newcommand\Directsum{\bigoplus}
\begin{enumerate}
\item $H^*_T(Y)$ is a free module over $H^*_T(pt)$.
\item $Y$ has a cell decomposition by complex cells $X_f$ (see
  \cite{BB}), corresponding to the fixed points $f\in Y^T$, whose
  closures give a basis of equivariant cohomology as a module over
  $H^*_T(pt)$.  The restriction $[\overline{X_f}]|_f$ of a class
  $[\overline{X_f}]$ to the point $f$ is the product of the weights in
  the normal bundle at the point $f$.
\item The forgetful map 
  $H^*_T(Y) \tensor_{H^*_T(pt)} \rationals \to H^*(Y)$
  is an isomorphism on rational cohomology.
\item The natural map $H^*_T(Y) \to H^*_T(Y^T) \iso \Directsum_{Y^T}
  H^*_T(pt)$ is {\em injective}\footnote{This is very odd from the
    point of view of ordinary cohomology -- we are restricting
    cohomology classes to individual points, which would seem very
    forgetful, but the above assertion says that we actually lose no
    information. The intuition should be that the {\em equivariant}
    cohomology of a point is very big, being a polynomial ring.}.
\end{enumerate}

These last two statements are the most combinatorially interesting:
they say that we can do all our calculations with lists of polynomials
indexed by the fixed points $Y^T$, and then once we understand the ring
$H^*_T(Y)$, we can (if we wish) set all the generators of $H^*_T(pt)$ 
to zero and recover ordinary cohomology. 

It remains to understand the image of $H^*_T(Y)$ inside
$\Directsum_{Y^T} H^*_T(pt)$.  
Given a class $\alpha \in H^*_T(Y)$ and a fixed point $f\in Y^T$, let
$\alpha|_f \in H^*_T(pt)$ denote the restriction of $\alpha$ to $f$,
so that this map is $\alpha \mapsto (\alpha|_f)_{f\in Y^T}$.
As it turns out, this image can be
characterized by the $T$-invariant copies of $\CP^1$ inside $Y$.

\begin{Theorem*}\cite{GKM}
  Let $T$ act on a smooth projective variety $Y$, with $Y^T$ finite,
  and let $Z$ be a $T$-invariant $\CP^1$ inside $Y$.  Then $Z^T$
  consists of two points $Z_{north}$ and $Z_{south}$.  If $w_Z$ is the
  weight of the $T$-action on the one-dimensional tangent space
  $T_{Z_{south}} Z$, then we have the ``GKM condition''
$$ \alpha|_{Z_{south}} - \alpha|_{Z_{north}} 
   \hbox{ is a multiple of } w_Z $$ 
for all classes $\alpha \in H^*_T(Y)$.

Conversely, suppose that there are only finitely many $T$-invariant
$\CP^1$'s, and $(\alpha|_\lambda)$ is an element of $\Directsum_{Y^T}
H^*_T(pt)$ which obeys the GKM condition for every $T$-invariant
$\CP^1$.  Then $(\alpha|_\lambda)$ lies in the image of $H^*_T(Y)$
inside $\Directsum_{Y^T} H^*_T(pt)$.
\end{Theorem*}

The first statement can be proven by applying the functor $H^*_T$ to
the inclusions $$\{Z_{north},Z_{south}\} \into Z \into Y.$$
The converse is deeper; see \cite{GKM}.

\subsection{Grassmannians.}\label{ssec:Gr}

We now apply this technology to the case of the Grassmannian $\Grkn$,
in order to verify our claimed combinatorial description of
$H^*_T(\Grkn)$.

The torus $T$ acting in this case is the $n$-dimensional torus $T =
U(1)^n$, so our base ring $H^*_T(pt)$ is $\integers[y_1,\ldots,y_n]$.
This torus acts by the diagonal action on $\complexes^n$, and thus
also acts on $\Grkn$.  The fixed points $\Grkn^T$ are just the
coordinate $k$-planes $\{ \complexes^\lambda: \lambda \in {n \choose k} \}$; 
we shall abuse notation and refer to the fixed point
$\complexes^\lambda$ simply as $\lambda$.

Two fixed points $\lambda$, $\lambda'$ are
connected by a $T$-invariant $\CP^1$ if and only if $\lambda = (i
\leftrightarrow j) \lambda'$ for some $1 \leq i < j \leq n$.  If $Z$
is such a $\CP^1$, then the action of $T$ on $T_{\lambda} Z$ 
has weight $\pm(y_j-y_i)$, and similarly for $T_{\lambda'} Z$.

From the GKM theorem we thus see that the equivariant cohomology
ring $H^*_T(\Grkn)$ is isomorphic (both as a ring and as a
$H^*_T(pt)$-module) to the ring $H^*_T(\Grkn)$ defined earlier
combinatorially, and we shall no longer bother to distinguish the two rings.

For each $\lambda \in {n \choose k}$, the Schubert cycle $X_\lambda$
defined in the introduction is oriented and $T$-invariant, so it
induces an equivariant cohomology class $\tS_\lambda := [X_\lambda]$.
These cycles are the closures of a cell decomposition of $\Grkn$ 
into complex cells.
At the fixed point $\lambda$, 
the weights of the $T$ action on the normal bundle of $X_\lambda$ are given by 
$\{y_j - y_i: (i,j) \in \inv(\lambda)\}$, so we have
$\tS_\lambda|_\lambda = \prod_{(i,j) \in \inv(\lambda)} (y_j - y_i)$.
The only fixed points in $X_\lambda$ are those corresponding to
strings above $\lambda$, so $\tS_\lambda$ is supported above
$\lambda$.  From degree considerations we see that $\deg(\tS_\lambda)
= l(\lambda)$ (treating the $y_i$ as having degree 1), and from the
first half of the GKM theorem we see that $\tS_\lambda$ obeys the GKM
conditions.  Combining this with lemma \ref{lem:uniq} we see that
$\tS_\lambda$ is indeed the Schubert class corresponding to $\lambda$.
Furthermore, we see that $\tS_\lambda$ maps to the ordinary cohomology
class $S_\lambda$ under the forgetful map to ordinary cohomology.

This concludes our construction of the equivariant Schubert 
classes $\tS_\lambda$. In section \ref{sec:ms} we will also relate 
these classes to factorial Schur functions, which are polynomials 
in many more variables\footnote{It is remarkable that these lists
of polynomials can be wrapped up into individual polynomials. 
This can also be traced to a geometrical fact, which is that
Grassmannians can be constructed as symplectic quotients of affine
space. Then the Kirwan map from equivariant cohomology of affine space
(a polynomial ring) maps onto the equivariant cohomology of the
Grassmannian. Since we will always work with classes, rather than
factorial Schur functions, we do not go into the details of this
argument.}.

\subsection{Schubert classes form a basis.}

Having used topological considerations to construct the equivariant
Schubert classes $\tS_\lambda$, we shall use a simple combinatorial
argument to show that they form a basis for $H^*_T(\Grkn)$.

\begin{Proposition}\label{prop:basis}
  The $\tS_\lambda$ form a $H^*_T(pt)$-basis for $H^*_T(\Grkn)$. More
  specifically, any class $\alpha \in H^*_T(\Grkn)$ can be written
  uniquely as an $H^*_T(pt)$-linear combination of $\tS_\lambda$ using
  only those $\lambda$ such that $\lambda \geq \mu$ for some $\mu \in
  \supp(\alpha)$.
\end{Proposition}

\begin{proof}
This essentially follows from the upper-triangularity of the
Schubert classes with respect to the order on $n\choose k$.
We now give the details.

First, we show the $\tS_\lambda$ are linearly independent.  Suppose
for contradiction that $\sum_\lambda Y_\lambda \tS_\lambda = 0$ for
some $Y_\lambda \in H^*_T(pt)$ which are not identically zero. Among
all $\mu \in {n \choose k}$ with $Y_\mu \neq 0$, we pick a $\mu$ which
is minimal in the lattice $n \choose k$.  But then the restriction of
$\sum_\lambda Y_\lambda \tS_\lambda$ to $\mu$ is $Y_\mu\, \tS_\mu|_\mu
\neq 0$, contradiction.

To see that the $\tS_\lambda$ span, let $\alpha$ be a class one is
attempting to write as a $H^*_T(pt)$-linear combination of some
classes $\{\tS_\lambda\}$ satisfying those conditions. Let $\mu$ be a
minimal element of the support of $\alpha$. From the GKM conditions we
see that $\alpha|_\mu$ must be a multiple $\beta$ of $\prod_{(i,j) \in
  \inv(\mu)} (y_j - y_i) = \tS_\mu|_\mu$. Subtracting $\beta \tS_\mu$,
we can inductively reduce the support of $\alpha$ upwards until it is
gone. This only uses those $\tS_\lambda$ for which $\lambda \geq \mu$
for some $\mu \in \supp(\alpha)$.
\end{proof}

\begin{figure}[htbp]
  \begin{center}
    \epsfig{file=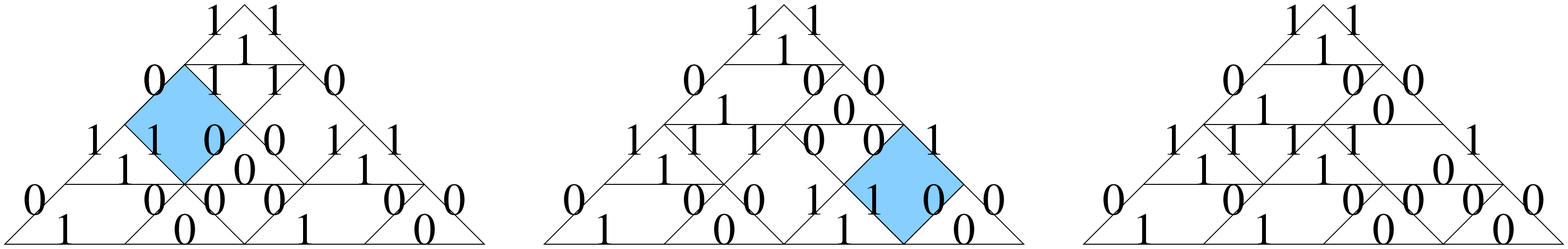,width=3.4in, height=1in}
    \caption{The puzzles computing
      $\tS_{0101} \tS_{1010} = (y_4-y_1) \tS_{1010} + \tS_{1100}$.
      Note that $y_4-y_1$ appears as $(y_3-y_1)+(y_4-y_3)$.
      Dualizing these puzzles, we get the three that compute 
      $\tS_{1010}\tS_{0101} = ((y_4-y_2)+(y_2-y_1)) \tS_{1010} + \tS_{1100}$.}
    \label{fig:S0101S1010}
  \end{center}
\end{figure}

\Ex Consider the class $\tS_{0101} \tS_{1010}$,
which is supported above $1010$ (refer back to figure \ref{fig:Gr24}
to see these classes). Following the algorithm given in 
proposition \ref{prop:basis} to write this in the Schubert basis,
we first subtract off a multiple of $\tS_{1010}$ itself, the
multiple being $\tS_{0101}|_{1010} = y_4-y_1$. 
The remainder is supported at $1100$, and is in fact $1$ times $\tS_{1100}$. 
In all
$\tS_{0101} \tS_{1010} = (y_4-y_1) \tS_{1010} + \tS_{1100}$.
This is an example of the ``equivariant Pieri rule'' proved in 
proposition \ref{eqvt-pieri}.
We verify the puzzle rule in this example, in figure \ref{fig:S0101S1010}.

Since the $\tS_\lambda$ form a $H^*_T(pt)$-basis of the ring
$H^*_T(\Grkn)$, we can define structure constants $\clmn \in
H^*_T(pt)$ for all $\lambda, \mu, \nu \in {n \choose k}$ by the
formula
$$ \tS_\lambda \tS_\mu = \sum_\nu \clmn \tS_\nu.$$
We record some elementary properties of these structure constants.

\begin{Lemma}\label{structure-lem}  
  The structure constant $\clmn$ has degree
  $l(\lambda)+l(\mu)-l(\nu)$, and vanishes unless $\nu \geq \lambda,
  \mu$ and $l(\nu) \leq l(\lambda) + l(\mu)$.  In the special case
  $l(\nu) = l(\lambda) + l(\mu)$, the structure constants are integers
  and agree with those from ordinary Schubert calculus.
  
  At the other extreme, when $\lambda = \nu$ 
  we have $c_{\lambda\mu}^\lambda = \tS_\mu|_\lambda$.
\end{Lemma}

\begin{proof}
  The first claim follows since each $\tS_\lambda$ has degree
  $l(\lambda)$.  In particular $\clmn$ vanishes when $l(\nu) >
  l(\lambda) + l(\mu)$.  The class $\tS_\lambda \tS_\mu$ is supported
  above $\lambda$ and above $\mu$, so by proposition \ref{prop:basis}
  we obtain the conditions $\nu \geq \lambda$.  If we apply the
  forgetful map from equivariant cohomology to ordinary cohomology,
  then the structure constants of non-zero degree all vanish, leaving
  only those with $l(\nu) = l(\lambda) + l(\mu)$, which explains the
  last claim in the first paragraph.
  
  Finally, if we restrict $\tS_\lambda \tS_\mu$ to $\lambda$, we
  obtain
$$ \tS_\lambda|_\lambda \tS_\mu|_\lambda 
= c_{\lambda \mu}^\lambda \tS_\lambda|_\lambda 
+ \sum_{\nu: \nu \neq \lambda} \clmn \tS_\nu|_\lambda.$$
Since $\clmn$ vanishes unless $\nu \geq \lambda$, and $\tS_\nu$ is
supported above $\nu$, we see that all the terms in the summation
vanish.  The claim then follows since $\tS_\lambda|_\lambda =
\prod_{(i,j) \in \inv(\lambda)} (y_j - y_i)$ is non-zero.
\end{proof}

From the above lemma we see that the equivariant structure constants
$\clmn$ compute Schubert classes when $l(\nu)=l(\lambda)$, 
and ordinary structure constants when $l(\nu)=l(\lambda)+l(\mu)$.  In the
next section we prove a Pieri rule which bridges the gap between these
two extreme cases.

\section{Pieri-based recurrence relations}\label{sec:pieri}

In this section we assume that $0 < k < n$, since the Schubert
calculus for the $k=0$ and $k=n$ cases are trivial.

Let $\dv := 000\ldots 010111\ldots 1$ denote the unique element of 
$n \choose k$ with one inversion. The corresponding Schubert class 
$\tS_\dv$ is the only one of degree $1$, coming from the unique 
Schubert divisor (hence the name). With $\tS_\dv$, and the associativity 
of the equivariant cohomology ring, we will establish recurrence
relations on the Schubert classes $\{\tS_\lambda\}$ and the
equivariant structure constants $c_{\lambda\mu}^\nu$.

\begin{Lemma}\label{lem:div}
  The Schubert divisor class $\tS_\dv$ is given by
  $$ \tS_\dv|_\lambda := \sum_{j=1}^n \lambda_j y_j - \sum_{i=1}^n y_i. $$
\end{Lemma}

\begin{proof}
  The right-hand side is clearly homogeneous of degree 
  $1 = \deg(\dv)$, supported above $\dv$, 
  and equals $y_{k+1} - y_k = \prod_{(i,j) \in \inv(\dv)} (y_j - y_i)$
  when restricted to $\lambda = \dv$.  It can easily be shown to also
  obey the GKM conditions.  The claim then follows from lemma
  \ref{lem:uniq}.
\end{proof}

Write $\lambda' \covers \lambda$ if $\lambda' > \lambda$ and
$l(\lambda') = l(\lambda) + 1$; this is the covering relation in the
lattice $n\choose k$.  Equivalently, $\lambda' \covers \lambda$ if
$\lambda' = \alpha 10 \beta$ and $\lambda = \alpha 01 \beta$ for some
strings $\alpha$, $\beta$.  
Thus for instance $110101 \covers 101101,110011$.

\begin{Proposition}[The equivariant Pieri rule]\label{eqvt-pieri}
  $$ \tS_\dv \tS_\lambda = (\tS_\dv|_\lambda) \tS_\lambda 
        + \sum_{\lambda': \lambda' \covers \lambda} \tS_{\lambda'}. $$
\end{Proposition}

\begin{proof}
From lemma \ref{structure-lem} and the fact that $\deg(\dv) = 1$ we have
$$\tS_\dv \tS_\lambda = (\tS_\dv|_\lambda) \tS_\lambda 
     + \sum_{\lambda': \lambda' \covers \lambda} 
     c_{\dv,\,\lambda}^{\lambda'} \tS_{\lambda'},$$
  where the $c_{\dv,\,\lambda}^{\lambda'}$ are the structure
  constants for ordinary Schubert calculus.  The claim then follows
  from the ordinary-cohomology Pieri rule\footnote{The ``Pieri rule''
    sometimes refers to a more general rule than we need here, for
    multiplying by $S_\lambda$ where $\lambda = 0\ldots 010\ldots
    01\ldots 1$. The equivariant version of this rule was recently
    formulated \cite{G}-positively in \cite{R}, for flag manifolds
    (not just Grassmannians).}  
  $ S_\dv S_\lambda = \sum_{\lambda': \lambda' \covers \lambda} S_{\lambda'}$
  (as proved in \cite{F}).
\end{proof}

In the appendix we shall give an alternative proof of proposition
\ref{eqvt-pieri} which does not go through the ordinary Pieri rule.

The equivariant Pieri rule gives a recurrence relation on
the structure constants $c_{\lambda \mu}^\nu$:

\begin{Theorem}\label{thm:crecurrence}\cite{MS}
For any $\lambda, \mu, \nu$ we have the recurrence relation 
$$ (\tS_\dv|_\nu - \tS_\dv|_\lambda) c_{\lambda\mu}^\nu = 
\bigg( \sum_{\lambda': \lambda'\covers \lambda} \ \ c_{\lambda'\mu}^\nu
- \sum_{\nu': \nu\covers\nu'} c_{\lambda\mu}^{\nu'} \bigg) $$
\end{Theorem}

The above recurrence was proven in \cite{MS} by a different argument;
it had also been observed by A. Okounkov.

\begin{proof}
We use associativity of multiplication in $H^*_T(\Grkn)$ and the equivariant 
Pieri rule to expand $\tS_\dv \tS_\lambda \tS_\mu$ in two different ways:
$$ (\tS_\dv \tS_\lambda) \tS_\mu 
= \big((\tS_\dv|_\lambda) \tS_\lambda 
        + \sum_{\lambda': \lambda'\covers\lambda} \tS_{\lambda'}\big)  \tS_\mu
= (\tS_\dv|_\lambda) \sum_\rho c_{\lambda\mu}^\rho \tS_\rho
  + \sum_{\lambda': \lambda'\covers\lambda} \sum_\rho c_{\lambda'\mu}^\rho \tS_\rho $$
and
$$ \tS_\dv (\tS_\lambda \tS_\mu) 
= \tS_\dv \sum_\rho c_{\lambda\mu}^\rho \tS_\rho
= \sum_\rho c_{\lambda\mu}^\rho \big( (\tS_\dv|_\rho) \tS_\rho 
        + \sum_{\rho: \rho'\covers\rho} \tS_{\rho'} \big).$$
Comparing coefficients of $\tS_\nu$, we get
$$ \tS_\dv|_\lambda\ c_{\lambda\mu}^\nu 
 + \sum_{\lambda': \lambda'\covers \lambda} c_{\lambda'\mu}^\nu 
 = c_{\lambda\mu}^\nu\, \tS_\dv|_\nu 
 + \sum_{\nu': \nu',\nu\covers\nu'} c_{\lambda\mu}^{\nu'} $$
as desired.
\end{proof}

The above recurrence gives us a purely combinatorial way to verify
that a putative formula for equivariant structure constants indeed works:

\begin{Corollary}\label{cor:crecurrence}  
  Let $0 \leq k \leq n$. Suppose that we have an assignment 
  $(\lambda, \mu, \nu) \mapsto d_{\lambda \mu}^\nu$ 
  from ${n\choose k}^3$ to $H^*_T(pt)$ obeying the following identities:
\begin{itemize}
\item  For any $\lambda \in {n \choose k}$, we have 
\begin{equation}\label{c1}
  d_{\lambda \lambda}^\lambda = \prod_{(i,j) \in \inv(\lambda)} (y_j - y_i).
\end{equation}

\item  For any $\lambda, \mu \in {n \choose k}$, we have
\begin{equation}\label{c2}
(\tS_\dv|_\lambda - \tS_\dv|_\mu) d_{\lambda \mu}^\lambda  = 
        \sum_{\mu': \mu'\covers\mu} d_{\lambda \mu'}^\lambda.
\end{equation}

\item  For any $\lambda, \mu, \nu \in {n \choose k}$ we have 
\begin{equation}\label{c3}
(\tS_\dv|_\nu - \tS_\dv|_\lambda) d_{\lambda\mu}^\nu = 
\sum_{\lambda': \lambda'\covers \lambda} d_{\lambda'\mu}^\nu
- \sum_{\nu': \nu\covers\nu'} d_{\lambda\mu}^{\nu'}.
\end{equation}
\end{itemize}
Then $\clmn = d_{\lambda \mu}^\nu$ for all $\lambda$, $\mu$, $\nu$.
\end{Corollary}

The identity \eqref{c1} thus
involves only $\Delta_{\lambda \lambda}^\lambda$-puzzles, while \eqref{c2}
involves $\Delta_{\lambda \mu}^\lambda$-puzzles and \eqref{c3}
involves general $\Delta_{\lambda \mu}^\nu$-puzzles.

\begin{proof}  
  If $k = 0$ or $k = n$ then we must have $\lambda = \mu = \nu$, and
  the claim follows from \eqref{c1}, Lemma \ref{structure-lem}, and
  the definition of $\tS_\lambda$.  So we assume $0 < k < n$.
  
  To begin with, we use the first two properties of $d$ to show that
  $d_{\lambda \mu}^\lambda = \tS_\mu|_\lambda$ (which we already
  knew to be equal to $c_{\lambda \mu}^\lambda$).  We induct on the
  quantity $l(\lambda) - l(\mu)$, which is clearly bounded from below.
  If $\lambda \neq \mu$, then $\tS_\dv|_\lambda - \tS_\dv|_\mu$ is
  non-zero, and the claim follows from \eqref{c2}, proposition
  \ref{eqvt-pieri}, and the induction hypothesis\footnote{Observe that
    this induction argument implies that $d_{\lambda \mu}^\lambda$
    vanishes unless $\lambda \geq \mu$, which is of course consistent
    with the support properties of $\tS_\mu$.  Similarly, the argument
    in the next paragraph shows that $d_{\lambda \mu}^\nu$ vanishes
    unless $\nu \geq \lambda$.}  (observing that
  $l(\lambda) - l(\mu') = l(\lambda) - l(\mu) - 1$).  In the
  base case $\lambda = \mu$ we instead use \eqref{c1} and our
  definition of the class $\tS_\mu$.
  
  Now we show that $\clmn = d_{\lambda \mu}^\nu$ in general.  We induct on the
  quantity $l(\nu) - l(\lambda)$, which is also bounded from below.
  The base case $\nu = \lambda$ follows from lemma
  \ref{structure-lem} and the previous paragraph.  In all other cases
  $\tS_\dv|_\nu - \tS_\dv|_\lambda$ is non-zero, and we can use
  \eqref{c3}, theorem \ref{thm:crecurrence}, and the observation that
  $l(\nu')-l(\lambda) = l(\nu)-l(\lambda') = l(\nu)-l(\lambda)-1$.
\end{proof}

We can tighten this further using the duality operation 
$P \mapsto P^*$ on puzzles, which takes $\partial P = \Delta_{\lambda \mu}^\nu$
to $\partial P^* = \Delta_{\mu^* \lambda^*}^{\nu^*}$.
Also, we have $wt(P^*) = \overline{wt(P)}$, where $x \mapsto \overline{x}$ 
is the involution on $H^*_T(pt)$ defined by
$\overline{y_i} := -y_{n+1-i}$ for $i = 1, \ldots, n$.  From these
observations and the definition of $d_{\lambda \mu}^\nu$ we see that
$$ d_{\lambda \mu}^{\nu} = \overline{d_{\mu^* \lambda^*}^{\nu^*}}.$$

\begin{Lemma}\label{lem:just1and3}
Let $\{ d_{\lambda \mu}^{\nu} \in H^*_T(pt) \}$ be a family satisfying
\eqref{c1}, \eqref{c3} and  
$ d_{\lambda \mu}^{\nu} = \overline{d_{\mu^* \lambda^*}^{\nu^*}}$.
Also assume the vanishing condition
\begin{equation}
  \label{vanishing}
  d_{\lambda\mu}^\nu = 0 \qquad \hbox{ unless $\nu\geq\lambda,\mu$.}
\end{equation}
Then \eqref{c2} follows automatically (and so by the corollary, 
 $\clmn \equiv d_{\lambda \mu}^\nu$).
\end{Lemma}

\begin{proof}
If we thus apply \eqref{c3} with $\lambda,\mu,\nu$ replaced by
$\mu^*,\lambda^*,\nu^*$ and apply the involution, we obtain
$$ (\overline{\tS_\dv|_{\nu^*}} - \overline{\tS_\dv|_{\lambda^*}})
d_{\lambda\mu}^\nu = \sum_{\mu': (\mu')^* \covers \mu^*}
d_{\lambda\mu'}^\nu - \sum_{\nu': \nu^*\covers (\nu')^*}
d_{\lambda\mu}^{\nu'}.
$$ 
From lemma \ref{lem:div} we have
$$ \overline{\tS_\dv|_{\nu^*}} - \overline{\tS_\dv|_{\mu^*}} =
\tS_\dv|_{\nu} - \tS_\dv|_{\lambda}$$ while from the definitions we
see that $(\mu')^* \covers \mu^*$ is equivalent to $\mu' \covers \mu$.
We obtain\footnote{This could also have been obtained from
\eqref{c3} and the commutativity property $d_{\lambda \mu}^\nu =
d_{\mu \lambda}^\nu$.  This commutativity property is of course true
(since it manifestly holds for the $c_{\lambda \mu}^\nu$) but is
non-trivial to prove using puzzles, so we rely instead on
duality.}
$$ (\tS_\dv|_{\nu} - \tS_\dv|_{\mu}) d_{\lambda\mu}^\nu =
\sum_{\mu': \mu' \covers \mu} d_{\lambda\mu'}^\nu - \sum_{\nu': \nu
\covers \nu'} d_{\lambda\mu}^{\nu'}.
$$
We now specialize this to the case $\nu = \lambda$, and use \eqref{vanishing}
to see that $d_{\lambda\mu}^{\lambda'} = 0$ when $\lambda \covers \lambda'$.  
The claim \eqref{c2} follows.
\end{proof}

We can now outline the proof of theorem \ref{thm:eqvtpuzcount}.  We
will assume $0 < k < n$ as the $k=0$, $k=n$ cases are trivial. To
prove the first conclusion of theorem \ref{thm:eqvtpuzcount}, it will
suffice by the above corollary to show that the quantity
$$ d_{\lambda \mu}^\nu 
        := \sum_{P: \partial P = \Delta_{\lambda \mu}^\nu} wt(P)$$
obeys the identities \eqref{c1}, \eqref{c3}, and \eqref{vanishing}.   
These will be proven in the next two sections.

\section{SW-NE rhombi, 
  and the proofs of \eqref{c1} and \eqref{vanishing}}\label{diag-sec}

In this section we prove the identities \eqref{c1} and \eqref{vanishing}.

\newcommand\flux{{flux}}
\newcommand\disc{{disc}}

We give first a ``Green's theorem'' argument to constrain the interior
of a puzzle from its boundary.  Suppose that $p$ is a SW-NE rhombus.
If we drop lines SE from $p$, they will poke out of the $j$th and
$(j+1)$st place of the South side of the puzzle for some $1 \leq j < n$.  
We then define the {\dfn discrepancy} of $p$ to be $\disc(p) := y_{j+1} - y_j$.

\begin{Lemma}\label{lem:Greendisc}
  Let $P$ be a $\Delta_{\lambda\mu}^\nu$-puzzle.  Then
$$ \sum_{\hbox{$p$ is a SW-NE rhombus of $P$}} \disc(p) 
        = \tS_\dv|_\nu - \tS_\dv|_\lambda.$$
\end{Lemma}

Note that the edges of a SW-NE rhombus are parallel to the $\lambda$
and $\nu$ sides of the puzzle.

\begin{proof}
  Let $p$ be any puzzle piece of $P$, and let $e$ be an edge of $p$.
  We give the pair $(p,e)$ a ``flux'' $\flux(p,e)$ as follows.  If
  $e$ is a $0$-edge, or a NW-SE $1$-edge, we set $\flux(p,e) := 0$.
  Otherwise we drop a line SE from $e$ until it pokes out of the $j$th
  place on the South side, and set $\flux(p,e) := \pm y_j$, where the
  sign $\pm$ is positive if $e$ is on the SW, SE, or S side of $p$,
  and negative if $e$ is on the N, NW, or NE side.
  
  Now compute the total flux $\sum_e \flux(p,e)$ of a puzzle piece
  $p$.  By checking each case from Figure \ref{fig:puzpieces} (and the
  equivariant piece) in turn, we see that $p$ has total flux zero
  unless $p$ is a NW-SE rhombus, in which case the total flux is
  $\disc(p)$.

  Finally, add up the flux of all the puzzle pieces in $P$. 
  At each internal edge, the contributions from the two pieces containing
  that edge cancel one another. So the total flux reduces to a sum over the
  edges on the boundary $\partial P$ of $P$, which can be computed as
  $$ \sum_{i=1}^n \nu_i y_i - \sum_{i=1}^n \lambda_i y_i 
        = \tS_\dv|_\nu - \tS_\dv|_\lambda.$$
  Combining this with the previous paragraph we obtain the lemma.
\end{proof}

This has some very pleasant corollaries:  

\begin{Corollary}\label{cor:fixedk}
  Let $P$ be a $\Delta_{\lambda\mu}^\nu$-puzzle.  

  Then $\lambda,\mu,\nu$ all have the same number of $1$'s (they are
  elements of the same $n\choose k$).

  Also, the number of
  rhombi in $P$ with edges parallel to the $\lambda$ and $\nu$ sides
  is $l(\nu)-l(\lambda)$.  Similarly when $\lambda$ is replaced with
  $\mu$ throughout.
\end{Corollary}

\begin{proof}
  For the first, specialize at $y_i \equiv 1$. (This argument can be
  presented much more simply than we have done here!) For the second, 
  specialize at $y_i \equiv i$.
\end{proof}

Another consequence is

\begin{Corollary}\label{cor:vanishing}
  Let $P$ be a $\Delta_{\lambda\mu}^\nu$-puzzle.  Then:
\begin{itemize}
\item We must have $\nu \geq \lambda$ and $\nu \geq \mu$ in the
  partial order on $n\choose k$. (This is \eqref{vanishing}.)
\item If $\lambda = \nu$, there can be no SW-NE rhombi.  
\item If $\lambda = \mu = \nu$, the can be no SW-NE or NW-SE rhombi.
\end{itemize}
\end{Corollary}

\begin{proof}
  Since the discrepancies $\disc(p)$ are all positive in the sense of
  \cite{G}, we see from the previous lemma that $\tS_\dv|_\nu -
  \tS_\dv|_\lambda$ is non-negative.  But this is equivalent to $\nu
  \geq \lambda$.  Furthermore, if $\nu = \mu$, then there cannot be
  any SW-NE rhombi, since $\tS_\dv|_\nu - \tS_\dv|_\lambda$ would then
  be strictly positive, a contradiction.

  To obtain the corresponding statements concerning $\mu$, we replace
  $P$ by the dual puzzle $P^*$ defined in the introduction.
  (Alternatively, one can ``dualize'' the proof of lemma
  \ref{lem:Greendisc} by the appropriate reflection and swapping of
  0-edges and 1-edges.)
\end{proof}

We now prove  \eqref{c1}, in the form of

\begin{Proposition}\label{prop:unique}
  There exists a unique $\Delta_{\lambda \lambda}^\lambda$-puzzle $P$,
  and its weight is 
$\prod_{(i,j) \in \inv(\lambda)} \, (y_j -  y_i)$.
\end{Proposition}

\begin{proof}
  Define a ``diamond'' in a puzzle to be any of the following objects:
\begin{itemize}
\item A N-S rhombus piece;
\item An equivariant puzzle piece;
\item Two triangular puzzle pieces joined by an E-W edge.
\end{itemize}
Note that the NW label on a diamond matches that on the SE,
likewise the NE and SW labels match.

  Let $P$ be a $\Delta_{\lambda\lambda}^\lambda$-puzzle. By the third
  conclusion of Corollary \ref{cor:vanishing}, $P$ contains no SW-NE
  or NW-SE rhombi.  Thus we can cut $P$ along all NW-SE and NE-SW
  lines without slicing through any rhombi. Except for the triangles
  attached to the South side, the sliced-up $P$ falls into diamonds.

  We analyze $P$ starting from the bottom. First, attach the isolated
  triangles. Then in each trough, fill in the unique diamond that
  fits.  We give the example of $\lambda = 1001$.\\
  \centerline{\epsfig{file=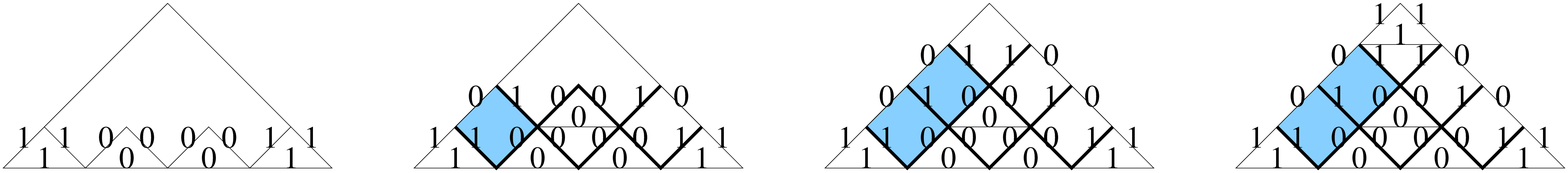,width=5in,height=1in}}
  Layer by layer, this creates the only puzzle with S edge $\lambda$
  (read left to right) that uses no NW-SE or NE-SW rhombi.  By the
  matching properties of diamonds, the NW and NE edges also end up
  labeled $\lambda$. This shows the existence and uniqueness.

  An equivariant piece comes whenever the trough to be filled has a
  $0$ on the SE and $1$ on the SW, coming from an inverted $0$ and $1$
  in $\lambda$. This shows that the weight is as advertised.
\end{proof}

It remains to prove \eqref{c3}.  This will be done in the next section,
at the end of which we give the proof of theorem \ref{thm:eqvtpuzcount}.

\section{Gashed puzzles: the proof of \eqref{c3}}\label{sec:cscab}

We give first the crucial definition, and then a rough indication
of the argument.

\begin{Definition}
  We define a {\dfn gashed puzzle} $(P,g)$ as a decomposition of a
  labeled equilateral triangle $\partial P$ into a collection $P$ of
  puzzle pieces, along with a line segment $g$ in the triangular lattice
  (which we refer to as the {\dfn gash}), such that
  \begin{itemize}
  \item The gash $g$ is contained in the equilateral triangle (either on
    the boundary $\partial P$ or in the interior), and is oriented either
    E-W or SW-NE;
  \item every edge {\em not} on the gash has at most one label (as in a
    non-gashed puzzle)
  \item if the gash is oriented SW-NE, then it is length $2$, and the
    labels on each side are a $0$ then a $1$ (read clockwise)
  \item if the gash is oriented E-W, then it is length at least $2$,
    with all but the first and last edge passing through the short diagonals 
    of some equivariant rhombi.  The labels on each side are a $0$,
    then the short diagonals of some equivariant rhombi, then $1$
    (read clockwise).
  \end{itemize}
\end{Definition}

\begin{figure}[htbp]
  \begin{center}
    \epsfig{file=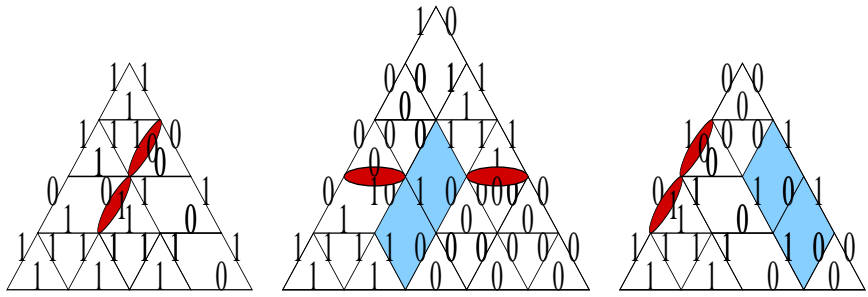,width=4.5in,height=1.5in}
    \caption{Three gashed puzzles. 
      The third has boundary $\Delta_{1010,0110}^{\ \ \ 1100}$, even
      though without the gash it would have been a puzzle with
      boundary $\Delta_{1100,0110}^{\ \ \ 1100}$.}
    \label{fig:gashex}
  \end{center}
\end{figure}

Some examples of gashed puzzles appear in figure \ref{fig:gashex}.  As
with non-gashed puzzles, we can define the {\dfn weight} $wt(P,g)$ of
a gashed puzzle to be the product of the weights of all the
equivariant pieces $p \in P$.  Thus for instance the second puzzle in
figure \ref{fig:gashex} has weight $(y_3-y_2)(y_4-y_2)$.

We now give an extremely rough indication of the argument, which decomposes 
\eqref{c3} into the four identities \eqref{left-ext}-\eqref{massage} to come.
Recall that equation \eqref{c3} computes $c_{\lambda\mu}^\nu$ from 
$\{ c_{\lambda'\mu}^\nu \}$ for $\lambda'\covers \lambda$ and 
$\{ c_{\lambda\mu}^{\nu'} \}$ for $\nu\covers \nu'$. We will take
puzzles with boundary $\Delta_{\lambda'\mu}^\nu$ and attach a gash on
their NW side, changing the boundary labels to $\Delta_{\lambda\mu}^\nu$.
(This will eventually give equation \eqref{left-ext}.)

\begin{figure}[htbp]
  \begin{center}
    \epsfig{file=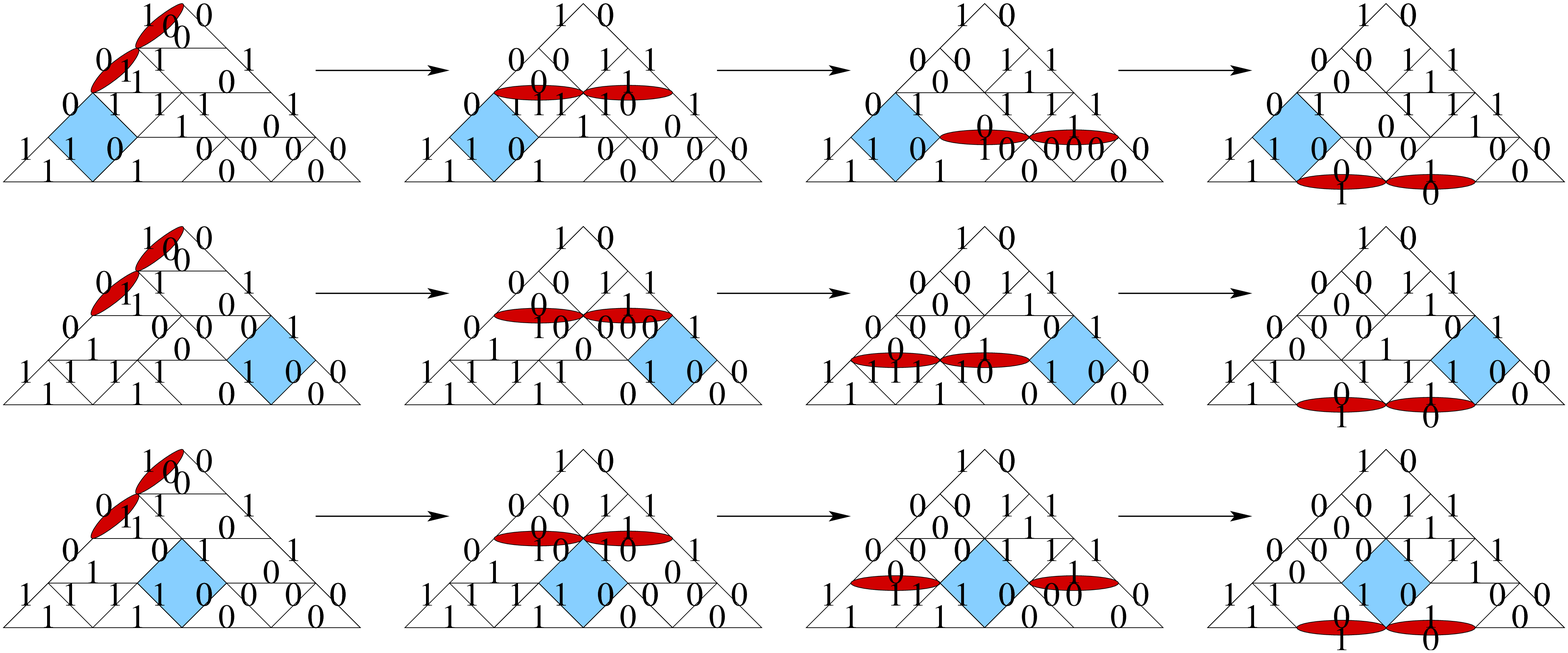,width=4.5in,height=3in}
    \caption{Gash propagation in gashed puzzles 
      with boundary $\Delta_{1010,0110}^{\,\,1100}$.
      In the third case the gash goes through an equivariant piece.
      Since there are no other gashed puzzles with this boundary,
      the left-hand side of \eqref{c3} vanishes, 
      and $c_{1010,0110}^{\ \ \ 1100}=0$.}
    \label{fig:crecEx1}
  \end{center}
\end{figure}

Then we will use some local rules for propagating a gash
through a gashed puzzle (the map $\phi$ in proposition \ref{prop:gashprop}, 
giving equation \eqref{massage}), preserving the weight.
The gash will usually come out on the S side, 
and when removed it leaves a puzzle with boundary $\Delta_{\lambda\mu}^{\nu'}$
(equation \eqref{right-ext}). If the gash {\em always} makes it through,
then 
$ \sum_{\lambda': \lambda'\covers \lambda} d_{\lambda'\mu}^\nu
- \sum_{\nu': \nu\covers\nu'} d_{\lambda\mu}^{\nu'} $
(the right-hand side of \eqref{c3}) will be zero.
This occurs in the $c_{1010,0110}^{\ \ \ 1100}$ example
given in figure \ref{fig:crecEx1}, and therefore 
$c_{1010,0110}^{\ \ \ 1100} = 0$.

\begin{figure}[htbp]
  \begin{center}
    \epsfig{file=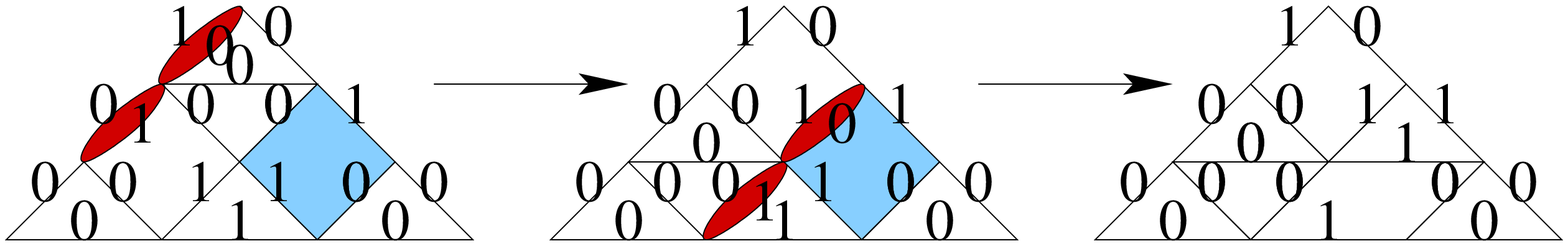,width=4in,height=1in}
    \caption{The gash propagates to an equivariant piece, 
      and they cancel one another, leaving a ``scab'' (defined in
      subsection \ref{ssec:scabs}).}
    \label{fig:crecEx2}
  \end{center}
\end{figure}

Frequently though, an equivariant piece can cause 
a gash to heal (or appear) on its own, and (by \eqref{telescope}) these
extra terms will give the left-hand side 
$(\tS_\dv|_\nu - \tS_\dv|_\lambda) d_{\lambda\mu}^\nu$
of equation \eqref{c3}.
The $c_{001,010}^{\ \ 010}=1$ example appears in figure \ref{fig:crecEx2}.

\begin{figure}[htbp]
  \begin{center}
    \epsfig{file=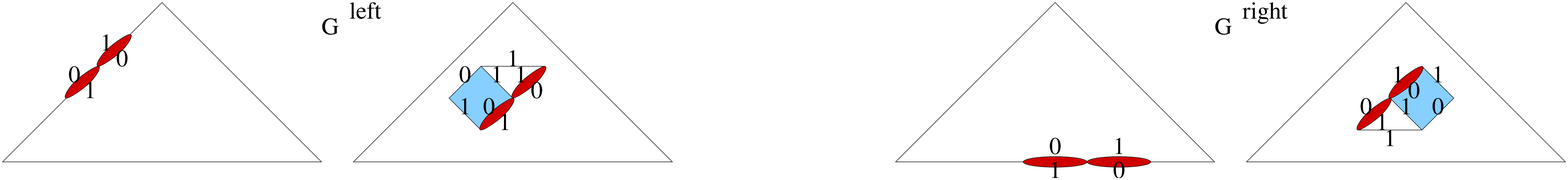,width=5in,height=1in}    
    \caption{Local pictures in gashed puzzles in 
      $G^{left}_{ext},G^{left}_{int},G^{right}_{ext},G^{right}_{int}$.}
    \label{fig:gashedpuz}
  \end{center}
\end{figure}

Fix $\lambda, \mu, \nu \in {n \choose k}$.  We define the set $G$ to
be the collection of all gashed puzzles $(P,g)$ with $\partial P =
\Delta_{\lambda \mu}^\nu$.  We define four subsets of $G$ (local pictures
given in figure \ref{fig:gashedpuz}):

\begin{align*}
 G^{left}_{ext} &:= \{(P,g) \in G : 
  \hbox{ $g$ lies in the NW boundary $P_{NW} = \lambda$ of $P$} \} &\hfill \\
 G^{left}_{int} &:= \{(P,g) \in G : 
  \hbox{ $g$ contains the SE edge of an equivariant piece} \} &\hfill \\
 G^{right}_{ext} &:= \{(P,g) \in G : 
  \hbox{ $g$ lies in the S boundary $P_{S} = w_0\cdot\nu$ of $P$} \} &\hfill \\
 G^{right}_{int} &:= \{(P,g) \in G : 
  \hbox{ $g$ contains the NW edge of an equivariant piece} \}.&\hfill
\end{align*}

We define $G^{left}$ to be the union of the (obviously disjoint) sets
$G^{left}_{int}$ and $G^{left}_{ext}$, and $G^{right}$ to be the union
of the (obviously disjoint) sets $G^{right}_{int}$ and
$G^{right}_{ext}$.  (The sets $G^{left}$ and $G^{right}$ may intersect.)

Our proof of \eqref{c3} will come down to four identities. Two are
very simple:

\begin{align}
  \sum_{(P,g) \in G^{left}_{ext}} wt(P,g) &= 
  \sum_{\lambda': \lambda' \covers \lambda} \, \sum_{P': \partial P' =
    \Delta_{\lambda' \mu}^\nu}  wt(P') \label{left-ext} \\
  \sum_{(P,g) \in G^{right}_{ext}} wt(P,g) &= 
  \sum_{\nu': \nu \covers \nu'} \, \sum_{P': \partial P' = \Delta_{\lambda
      \mu}^{\nu'}} wt(P') \label{right-ext} 
\end{align}

We begin with \eqref{left-ext}.
Let $(P,g)$ be an arbitrary element of $G^{left}_{ext}$.  Then $g$
lies on a pair of edges where $P_{NW}=\lambda$ reads 10.  If one
removes the gash, one obtains a non-gashed puzzle $P'$ with boundary
$\partial P' = \Delta_{\lambda' \mu}^\nu$, where $\lambda'$ is equal
to $\lambda$ but with 10 replaced by 01.  In particular, we have
$\lambda' \covers \lambda$.  This argument can be reversed; given any
$\Delta_{\lambda' \mu}^\nu$-puzzle $P'$ with $\lambda' \covers
\lambda$, we can take the two edges where $\lambda$ and $\lambda'$
disagree, and swap them to create a gashed puzzle $(P,g)$.  Since this
affects no equivariant pieces, one has $wt(P,g) = wt(P')$, 
and \eqref{left-ext} follows.

The proof of \eqref{right-ext} is exactly the same, except for one 
minor subtlety:
observe that if $g \subseteq P_S$ then the gash $g$ must have length 2, 
since $P_S$ does not contain the short diagonal of any equivariant
rhombus pieces.

The other two identities
\begin{align}
  \sum_{(P,g) \in G^{right}_{int}} wt(P,g) - \sum_{(P,g) \in
    G^{left}_{int}} wt(P,g) &= 
  \sum_{P': \partial P' = \Delta_{\lambda \mu}^\nu} wt(P')
  (\tS_\dv|_\nu - \tS_\dv|_\lambda)  \label{telescope} \\ 
  \sum_{(P,g) \in G^{right}} wt(P,g) &= 
  \sum_{(P,g) \in G^{left}} wt(P,g) \label{massage}
\end{align}
are more subtle and will be proved in the next two subsections.

\subsection{Proof of \eqref{telescope}.}\label{ssec:scabs}

To prove \eqref{telescope} we need to introduce the notion of a {\dfn scab}.  
Let $P'$ be a $\Delta_{\lambda \mu}^\nu$ puzzle.  We define a 
{\dfn left-scab} of $P'$ to be any pair $\kappa$ of puzzle pieces in
$P'$ consisting of an SW-NE rhombus sitting atop a downward
$1$-triangle.  Similarly define a {\dfn right-scab} of $P'$ to be any
pair $\kappa$ of puzzle pieces in $P'$ consisting of an upward
$1$-triangle sitting atop a SW-NE rhombus. 

We define the weight $wt(\kappa)$ of the scab by $wt(\kappa) :=
wt(p)$, where $p$ is the unique equivariant piece which can fit inside
the region occupied by $\kappa$.

\begin{Lemma*}  We have
  $$
  \sum_{(P,g) \in G^{left}_{int}} wt(P,g) 
  = \sum_{P': \partial P' = \Delta_{\lambda \mu}^\nu} \,
  \sum_{\kappa: \,\kappa \hbox{ is a left-scab of } P'} wt(P') wt(\kappa)
  $$
  and
  $$
  \sum_{(P,g) \in G^{right}_{int}} wt(P,g) 
  = \sum_{P': \partial P' = \Delta_{\lambda \mu}^\nu} \,
  \sum_{\kappa: \,\kappa \hbox{ is a right-scab of } P'} wt(P') wt(\kappa).
  $$
\end{Lemma*}

\begin{proof}
Let $(P,g)$ be a gashed puzzle in $G^{left}_{int}$.  Then $g$ must be
a NE-SW line segment of length 2, whose SW edge is the SE edge of an
equivariant piece $p$.  From Figure \ref{fig:puzpieces} we thus see
that there must be a downward 1-triangle $t$ between $p$ and the NE
edge of $g$, as in figure \ref{fig:gashedpuz}.
Observe that if we replace $p$ and $t$ with a left-scab
$\kappa$, we obtain an ungashed puzzle $P'$ with $\partial P' =
\partial P = \Delta_{\lambda \mu}^\nu$ and $wt(P') wt(\kappa) = wt(P')
wt(p) = wt(P)$.  \\
\centerline{\epsfig{file=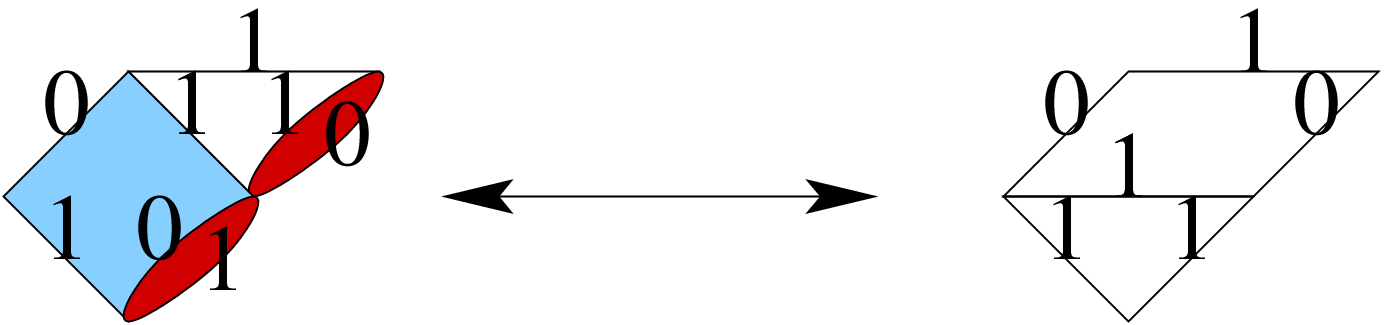,height=.7in,width=2in}}
This procedure can be reversed; given any ungashed
puzzle $P'$ with $\partial P = \Delta_{\lambda \mu}^\nu$, and given
any left-scab $\kappa$ of $P'$, we can replace the scab $\kappa$ with
an equivariant piece $p$ and a downward 1-triangle $t$, creating a
gashed puzzle $P$ with $\partial P = \partial P' = \Delta_{\lambda
\mu}^\nu$ and $wt(P) = wt(P') wt(p) = wt(P') wt(\kappa)$.  This proves
the first claim. 

The second claim is similar (indeed, it is essentially a $180^\circ$
rotation of the first claim) and is left to the reader.
\end{proof}

The argument above motivates the terminology; when a gash closes up it
leaves a scab, and conversely, a scab can come off producing a new gash.

From the above lemma, we see that to prove \eqref{telescope} it will
suffice to show

\begin{Lemma*}  Let $P'$ be a $\Delta_{\lambda \mu}^\nu$-puzzle.  Then
$$ \sum_{\kappa: \kappa \hbox{ is a right-scab of } P'} wt(\kappa) -
\sum_{\kappa: \kappa \hbox{ is a left-scab of } P'} wt(\kappa) =
\tS_\dv|_\lambda - \tS_\dv|_\nu.$$
\end{Lemma*}

\begin{proof}
This is another Green's theorem argument.
Let $p$ be any puzzle piece of $P'$, and let $e$ be an edge of $p$.
We give the pair $(p,e)$ a ``flux'' $\flux(p,e)$ as follows.  If
$e$ is a $0$-edge, or a NW-SE $1$-edge, we set $\flux(p,e) := 0$.  If
$e$ is an E-W $1$-edge, we drop a line SW from $e$ until it pokes out
of the $i$th place on the South side, and set $\flux(p,e) := +y_i$ if
$e$ is on the N side of $p$, and $\flux(p,e) := -y_i$ if $e$ is on the
S side of $p$.  Similarly, if $e$ is a SW-NE $1$-edge, we drop a line
SE from $e$ until it pokes out of the $j$th place on the South side,
and set $\flux(p,e) := +y_j$ if $e$ is on the NW side of $p$, and
$\flux(p,e) := -y_j$ if $e$ is on the SE side of $p$.

Now compute the total flux $\sum_e \flux(p,e)$ of a puzzle piece $p$.
By checking each case from Figure \ref{fig:puzpieces} (and the
equivariant piece) in turn, we see that $p$ has total flux zero unless
$p$ is a $1$-triangle.  Furthermore, if $p$ is an upward $1$-triangle
sitting atop the south boundary $P_S$, then $p$ also has total flux 0.
Finally, if $p$ is an upward $1$-triangle sitting atop a downward
$1$-triangle $p'$, then the total flux of $p$ and $p'$ is zero.

Thus the only upward $1$-triangles $p$ which have non-zero flux are
those which sit atop SW-NE rhombi.  But in that case $p$ belongs to a
right-scab $\kappa$, and the total flux of $p$ can be easily computed
to equal $wt(\kappa)$.  Similarly, only the downward $1$-triangles $p$
which have non-zero flux are those which sit below SW-NE rhombi, so
they belong to a left-scab $\kappa$, and the total flux of $p$ can be
easily computed to be $-wt(\kappa)$.

Finally, add up the flux of all the puzzle pieces in $P$.  At each
  internal edge, the contributions from the two pieces containing that
  edge cancel one another. So the total flux reduces to a sum over the
  edges on the boundary $\partial P$ of $P$, which can be computed as
$$ -\sum_{i=1}^n \nu_i y_i + \sum_{i=1}^n \lambda_i y_i =
\tS_\dv|_\lambda - \tS_\dv|_\nu.$$ Combining this with the previous
paragraph we obtain the lemma.
\end{proof}

\begin{proof}[Proof of \eqref{telescope}.]
This follows from the two lemmata just proven:
$$ 
  \sum_{(P,g) \in G^{right}_{int}} wt(P,g) - \sum_{(P,g) \in
    G^{left}_{int}} wt(P,g) 
  \qquad\qquad\qquad\qquad\qquad\qquad\qquad\qquad\qquad
$$
\begin{align*}
&=   \sum_{P': \partial P' = \Delta_{\lambda \mu}^\nu} wt(P')
   \bigg(\sum_{\kappa: \kappa \hbox{ is a left-scab of } P'} wt(\kappa) -
   \sum_{\kappa: \kappa \hbox{ is a right-scab of } P'} wt(\kappa)\bigg) \\
&=   \sum_{P': \partial P' = \Delta_{\lambda \mu}^\nu}
   wt(P')   (\tS_\dv|_\nu - \tS_\dv|_\lambda).
\end{align*}
\end{proof}

\subsection{Proof of \eqref{massage}.}

Equation \eqref{massage} is equivalent to
$$\sum_{(P,g) \in G \setminus G^{right}} wt(P,g) = \sum_{(P,g) \in G
\setminus G^{left}} wt(P,g).$$ This shall be an immediate consequence of

\begin{Proposition}  \label{prop:gashprop}
  There exists a weight-preserving bijection $\phi$ from $G \setminus
  G^{right}$ to $G \setminus G^{left}$.
\end{Proposition}

\begin{proof}
  Let $(P,g)$ be an element of $G \setminus G^{right}$.  
  We shall construct an element $(P',g') = \phi(P,g)$ of 
  $G \setminus G^{left}$ for which $wt(P',g') = wt(P,g)$.  
  This will require only a local surgery on $P$, in which 
  some pieces are replaced and the gash moves.\footnote{It is worth
    noting, for readers cognizant of the ``gentle path''
    technology of \cite{Hon2}, that the center of the gash always
    moves along a gentle path. Indeed, the loop-breathing in \cite{Hon2}
    can be interpreted as introducing a ``double gash'' crossing the gentle
    loop, propagating one gash around the loop, and once it gets back
    removing them both.}
  
  Suppose first that $g$ is a SW-NE gash, and consider the pieces
  to its right, with a vertex on the center of the gash.
  Since $(P,g) \not \in G^{right}_{int}$, there cannot be an
  equivariant piece immediately to the right of $g$, which
  leaves three possibilities:\\ \vskip.05in
\centerline{  \epsfig{file=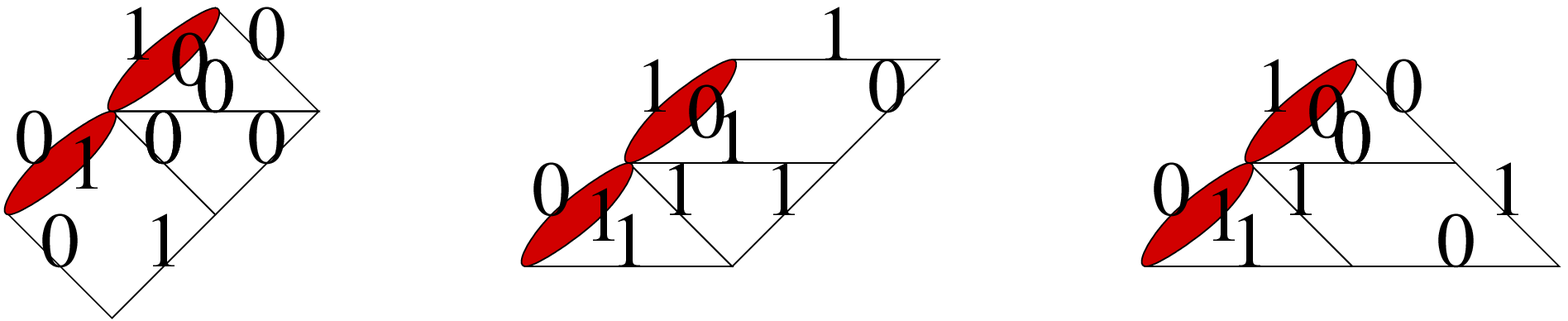,height=.7in,width=2.5in} }

  In each of these cases, we remove the pieces and gash,
  and replace them as follows:\\ \vskip.05in
\centerline{  \epsfig{file=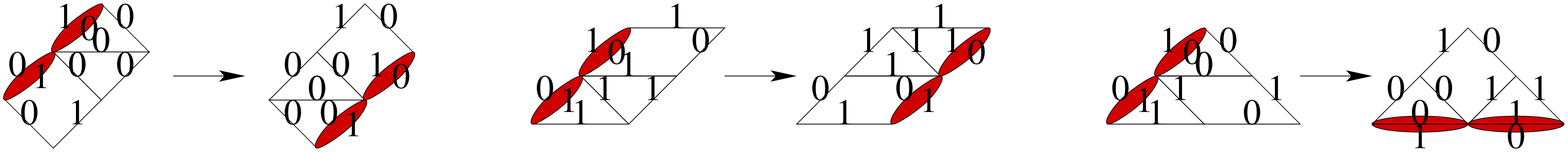,height=.7in,width=5in} }

In each case, the labels on the {\em boundary} do not change, which means
that the new set of pieces and gash match fit into the puzzle where the
old ones were. So this creates a new gashed puzzle $(P',g')$, 
and this is how we define $\phi(P,g)$.

Now take the case that $g$ is an E-W gash.  Since $(P,g) \not \in
G^{right}_{ext}$ we see that $g$ is not on the $S$ edge of the puzzle.

Suppose first that $g$ has length 2, and consider the pieces below $g$
with a vertex on the center of the gash. There are four possibilities,
which we give below, along with their replacements in $(P',g') =: \phi(P,g)$:
\\ \vskip.05in
\centerline{  \epsfig{file=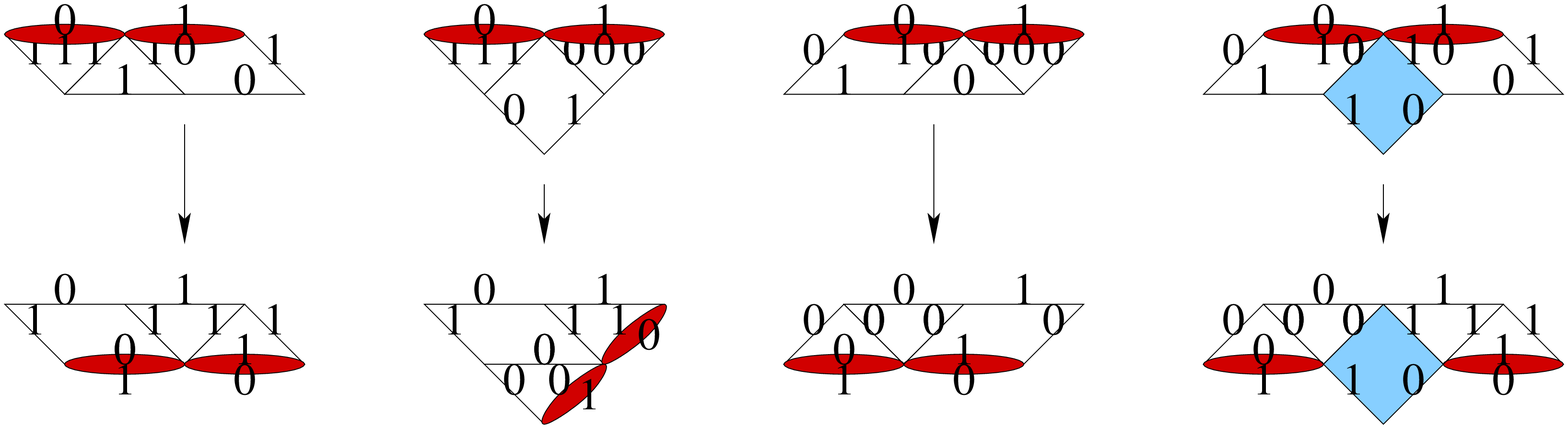,height=2in,width=4in} }

Now suppose that $g$ has length $l > 2$.  Then the two extreme edges of $g$ 
are separated by $l-2$ equivariant pieces as shown below (with $l=5$):
\\ \vskip.05in
\centerline{  \epsfig{file=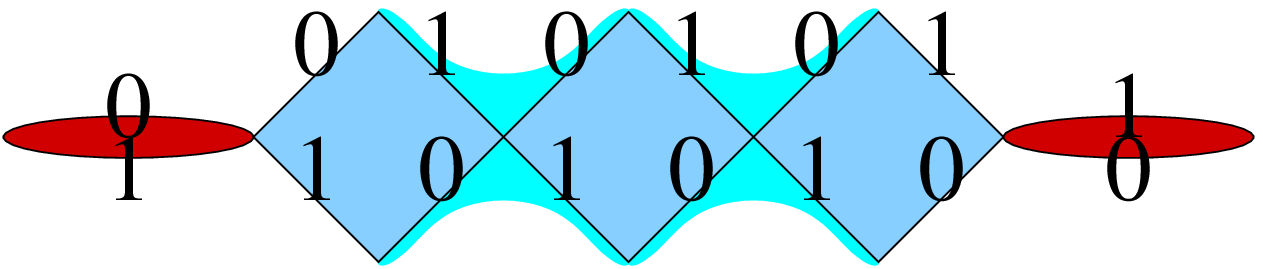,height=.7in,width=2in} }
Observe that between and below any two equivariant pieces on the gash
there must be another equivariant piece (since nothing else can fit in 
that space). That leaves $2^2 = 4$ possibilities for the ends, depending
on whether there are more equivariant pieces in that row below. 
In each case we move the
gash down one step, possibly stretching it or shrinking it by length 1:
\\ \vskip.05in
\centerline{  \epsfig{file=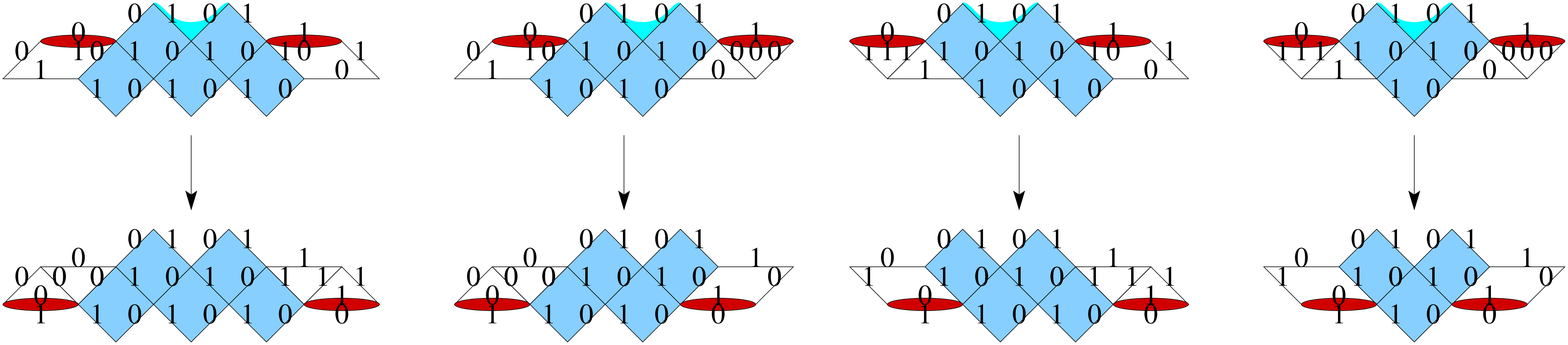,height=2in,width=6in} }
(In these pictures the gash begins length 4, and ends length 5, 4, 4, or 3.)
This creates a new gashed puzzle $(P',g')$, with which we define $\phi(P,g)$.

We have now defined $\phi(P,g) = (P',g')$ for all $(P,g) \in G
\setminus G^{left}$.  A simple examination of all cases verifies that
$(P',g')$ is a gashed puzzle with $\partial P' = \partial P =
\Delta_{\lambda \mu}^\nu$, and also $(P',g') \not \in G^{right}$.
It is obvious that $wt(P',g') = wt(P,g)$ since no equivariant pieces
are created, destroyed, or moved.

If we rotate these local-replacement recipes by $180^\circ$, we get
a similar map $\phi'$ from $G\setminus G^{right}$ to $G\setminus G^{left}$.
This is easily checked to be the inverse of $\phi$, which is therefore
a bijection.
\end{proof}

One can use the same rules to define a correspondence between
$G^{left}$ and $G^{right}$, but they must be iterated. This was the
viewpoint of the examples in figures \ref{fig:crecEx1} and \ref{fig:crecEx2}
at the beginning of the section.

\begin{proof}[Proof of theorem \ref{thm:eqvtpuzcount}.]
  Combining \eqref{left-ext}, \eqref{right-ext}, \eqref{telescope},
  \eqref{massage}, we obtain
$$
    \sum_{P': \partial P' = \Delta_{\lambda \mu}^\nu} wt(P')
    (\tS_\dv|_\nu - \tS_\dv|_\lambda)  
\qquad\qquad\qquad\qquad\qquad\qquad\qquad\qquad\qquad
\qquad\qquad\qquad\qquad
$$
  \begin{align*}
&= 
    \sum_{(P,g) \in G^{right}_{int}} wt(P,g) - \sum_{(P,g) \in
      G^{left}_{int}} wt(P,g) 
    &\hbox{by \eqref{telescope}} \\
 &=   \sum_{(P,g) \in G^{left}_{ext}} wt(P,g) - \sum_{(P,g) \in
    G^{right}_{ext}} wt(P,g) 
    &\hbox{by \eqref{massage}} \\
 &=   \sum_{\lambda': \lambda' \covers \lambda} \, \sum_{P': \partial P' =
    \Delta_{\lambda' \mu}^\nu}  wt(P') 
-  \sum_{\nu': \nu \covers \nu'} \, \sum_{P': \partial P' = \Delta_{\lambda
      \mu}^{\nu'}} wt(P') 
    &\hbox{by \eqref{left-ext} and \eqref{right-ext}}
  \end{align*}
  and this is \eqref{c3}.

  Then by lemma \ref{lem:just1and3} and proposition \ref{prop:unique},
  we obtain the first statement of theorem \ref{thm:eqvtpuzcount}.

  There is one foolish subtlety in obtaining the second statement
  of the theorem: the first statement 
  (and the recurrences \eqref{c2}-\eqref{c3})
  only constrain $d_{\lambda\mu}^\nu$ for $\lambda,\mu,\nu$ all having
  the same number of $1$s, and so a priori we might worry that the
  product $\tS_\lambda \tS_\mu$ might be miscalculated to have some
  extra terms in which $\nu$ has a different number of $1$s. But by
  corollary \ref{cor:fixedk} above, the number of $1$s is the same on
  all sides of a puzzle. The second statement follows.
\end{proof}

\section{The Molev-Sagan problem}\label{sec:ms}

In this section we compare the results of this paper with the earlier
work in \cite{MS}, which was a major source of inspiration for this
paper.  We also give a homological (or ``geometrical'') interpretation
of the structure constants computed in \cite{MS}.

The paper \cite{MS} is concerned with multiplying ``factorial Schur
functions'' $s_\lambda(x|y)$ for $\lambda \in {n \choose k}$.  These
functions are polynomials in two sets of variables $\{x_1, \ldots,
x_k\},\{y_1, \ldots y_n\}$, and are related to the classes
$\tS_\lambda$ by

\begin{Lemma*}\cite{O}, \cite{MS}  
For any $\lambda, \mu \in {n \choose k}$ 
we have $\tS_\lambda|_\mu = s_\lambda(y_\mu|y)$
where $y_\mu := \{ y_{i_1}, \ldots, y_{i_k} \}$ and $n \geq i_1 >
\ldots > i_k \geq 1$ are the $k$ integers $\{ 1 \leq i \leq n: \mu_i =
1 \}$ in decreasing order.
\end{Lemma*}

\begin{proof}
It is easy to check the GKM conditions, and these $s_\lambda$ have the
right vanishing and normalization conditions (as proved in \cite{O} and
repeated in \cite{MS}). Then apply lemma \ref{lem:uniq}.
\end{proof}

The problem solved in \cite{MS} is more general than the one we have
stated: they consider the mixed structure constants
$e_{\theta\mu}^{\nu}(y,z)$, which are polynomials in variables $y$ and
$z$, given by the product expansion
$$ s_{\theta}(x|z) \, s_{\mu}(x|y)= \sum_{\nu} \,
e_{\theta\mu}^{\nu}(y,z)\ s_{\nu}(x|y). $$ The $e_{\theta \mu}^\nu$
reduce to the structure constants for equivariant cohomology of
Grassmannians at the specialization $y \equiv z$ (and to ordinary
cohomology at $y\equiv z\equiv 0$).\footnote{%
In \cite{MS} they also permit $s_\theta$ to be a {\em skew} Schur function,
not just a Schur function, but we have not been able to find any
cohomological interpretation of these.}

The formula in \cite{MS} writes $e_{\theta\mu}^\nu(y,z)$ as a sum over
``barred tableaux,'' each one contributing a certain product $\prod
(y_i-z_j)$. In that sense their formula is positive (and reduces to
the Littlewood-Richardson rule, in the case that $l(\nu) =
l(\theta)+l(\mu)$). Unfortunately, many of their terms have
$(y_i-z_j)$ factors with $i\leq j$,
%
as the example
$$ s_{10}(x|z)\, s_{01}(x|y) = (x_1 - z_1)\cdot 1 = (x_1 - y_1) + (y_1
- z_1)\cdot 1 = s_{10}(x|y) + (y_1-z_1) s_{01}(x|y) $$ already shows.
%
For this reason, the computation of Molev-Sagan structure constants is
too general a setting for finding a formula (as in theorem
\ref{thm:eqvtpuzcount}) for equivariant Schubert calculus that is
manifestly positive in the sense of \cite{G}.

We now give a cohomological interpretation of the structure
constants $e_{\theta \mu}^\nu$, which we christen ``triple Schubert
calculus'', and sketch how one can also compute these co-efficients
using ``MS-puzzles''.  (Cohomology does not explicitly appear in
\cite{MS} -- they consider the computation of the $\{ e_{\theta \mu}^\nu\}$
purely as a combinatorial question.)

\subsection{Double Schubert calculus vs. equivariant Schubert calculus.}\label{double-eqvt-sec}

\newcommand\Fln{Flags(\complexes^n)} 
In this subsection we recall the (well-known) connection between
double Schubert calculus and equivariant Schubert calculus.  In a
nutshell, the connection is that $H^*_T(X) \onto H^*(\Fln \times X)$
for any partial flag manifold $X$, taking equivariant Schubert classes
to ``double Schubert classes''.
We begin by recalling some standard material on double Schubert
calculus (originally defined in \cite{LS}), and its geometric interpretation.

Let $\Fln$ denote the space of flags (i.e. maximal chains of
subspaces) in $\complexes^n$. This has a transitive action of $\GLn$
induced from its action on $\complexes^n$, and the stabilizer of the
standard flag $(F_i)$ is the upper triangular matrices $B$, so $\Fln
\iso \GLn/B$.

If we denote the lower triangular matrices by $B_-$, then the Schubert
cells $X_\lambda$ on the Grassmannian $\Grkn$ are exactly the $B_-$
orbits, whose Poincar\'e duals gave us the Schubert basis $S_\lambda$ of
ordinary cohomology $H^*(\Grkn)$. These were indexed by patterns 
$\lambda \in {n\choose k}$
recording the intersection of the $k$-plane with the anti-standard flag.

Analogously, we can consider closures of the $\GLn$-orbits on $\Fln
\times \Grkn$, which are again indexed by $n\choose k$ (recording the
intersection of the $k$-plane with the flag, which is now varying). The
Poincar\'e duals $D_\lambda \in H^*(\Fln \times \Grkn)$ of these
$\GLn$-orbit closures then form an $H^*(\Fln)$-basis for
$H^*(\Fln \times \Grkn)$ as $\lambda$ varies over $n \choose k$.

Since $H^*(\Fln \times \Grkn)$ is a ring as well as an
$H^*(\Fln)$-module, we can define structure constants
$\{f_{\lambda\mu}^\nu\} \in H^*(\Fln)$ for the multiplication:
$$ D_\lambda D_\mu = \sum_\nu f_{\lambda\mu}^\nu D_\nu.$$
The computation of the $f_{\lambda \mu}^\nu$ is the concern of {\dfn
double Schubert calculus}, and has the following homological
interpretation.  Fix a generic element $g$ of $\GLn$.  The class
$D_\lambda D_\mu$ corresponds to an irreducible cycle of pairs $(F,V)
\in \Fln \times \Grkn$ where $V$ satisfies two intersection conditions
with $F$: $V$ intersects $F$ $\lambda$-much, and intersects $gF$
$\mu$-much.  The class $\sum_\nu f_{\lambda \mu}^\nu D_\nu$
corresponds to a union of cycles, each of which put only one condition
on $V$ (that it intersect $F$ $\nu$-much), while also requiring that
$F$ live in a cycle Poincar\'e dual to $f_{\lambda\mu}^\nu$. So the
equation requiring these two to be homologous is somehow splitting the
double burden on $V$ to a single burden on $V$ and a single burden on
$F$.

Restricting $\Fln \times \Grkn$ to $\Fln \times pt$, each $D_\lambda$
maps to $S_\lambda$, which shows that these $f_{\lambda\mu}^\nu$
generalize the structure constants $c_{\lambda \mu}^\nu$ of ordinary
Schubert calculus.  This was also true of the structure constants of
equivariant Schubert calculus, and like them, the $f_{\lambda \mu}^\nu$
carry a degree.\footnote{However, this degree is the degree of a
cohomology class rather than a polynomial, and to be precise it is
$2(l(\lambda)+l(\mu)-l(\nu))$ rather than $l(\lambda)+l(\mu)-l(\nu)$.}

\newcommand\dom{\backslash} 
We can connect double Schubert calculus
with equivariant Schubert calculus using the following property
of equivariant cohomology: if $G\times K$ acts on $X$, and $K$'s
action is free, then $H^*_{G\times K}(X) \iso H^*_G(X/K)$.  Letting
$P$ denote the stabilizer of the standard $k$-plane $\complexes^{w_0\cdot
\id}$, so $\Grkn \iso \GLn/P$, we thus have
$$ H^*_T(\Grkn) \iso H^*_{T\times P}(\GLn) $$
$$ \iso H^*_{T\times\GLn\times P}(\GLn\times \GLn) \iso H^*_\GLn(T\dom
\GLn \times \GLn/P) $$
$$ \onto H^*(T\dom \GLn \times \GLn/P) \iso H^*(\Fln \times \Grkn),$$
where we have used the fact that $\GLn$ is isomorphic to $\GLn \times
\GLn$ quotiented by the diagonal action of $\GLn$.  Thus the structure
constants $c_{\lambda \mu}^\nu$ for equivariant Schubert calculus live in
$H^*_\GLn(\Fln) \iso H^*_{\GLn\times B}(\GLn) \iso H^*_B(pt) \iso H^*_T(pt)$, 
and the forgetful map from $H^*_\GLn(\Fln)$ to $H^*(\Fln)$
maps these constants to the double Schubert calculus constants
$f_{\lambda \mu}^\nu$.

To summarize the above discussion, while equivariant Schubert calculus
essentially lacks a definition in terms of intersecting cycles, 
one is provided by double Schubert calculus, of which equivariant
Schubert calculus is a refinement.\footnote{One can show that it
  is the only such refinement satisfying some natural stability
  properties in the limit $n \to \infty$. This is perhaps a bad
  way to see things, though, since such a limit can only be defined
  for classical Lie groups, whereas double and equivariant
  Schubert calculus can be defined for arbitary Lie groups.}

\subsection{Triple Schubert calculus.}

Consider the cohomology ring $H^*\big(\Fln \times \Grkn \times
\Fln\big)$ as a module over the cohomology ring $H^*(\Fln\times \Fln)$ of the
first and third factors.  Since the classes $D_w$ form a basis of $H^*(\Fln
\times \Grkn)$, we see that the classes $\{D_w \tensor 1\}$ form a
$H^*(\Fln\times\Fln)$-basis for $H^*\big(\Fln \times \Grkn \times \Fln\big)$.
One could compute the structure constants for multiplication in this
basis, but one would just obtain the double Schubert constants 
$f_{\lambda \mu}^\nu$ again (or to be pedantic, we would obtain
 $f_{\lambda\mu}^\nu \tensor 1$).
 
 Since $H^*(\Fln \times \Grkn) \equiv H^*(\Grkn \times \Fln)$, the
 classes $D_w$ induce a corresponding basis $D_w^T$ for $H^*(\Grkn
 \times \Fln)$.  The classes $1 \tensor D_w^T$ thus form another
 $H^*(\Fln\times\Fln)$-basis for $H^*\big(\Fln \times \Grkn \times
 \Fln\big)$.  Again, the structure constants for this basis are no
 richer than for double Schubert calculus.

In \cite{MS} the authors (implicitly) considered the hybrid problem of
computing the structure constants ${e''}_{\theta \mu}^\nu \in H^*(\Fln
\times \Fln)$ in the expansion
$$ (1 \tensor D_\theta^T) (D_\mu \tensor 1) = \sum_\nu
{e''}_{\theta\mu}^\nu (D_\nu \tensor 1).$$ This has the following
homological interpretation. We are now looking for triples $(F_1,V,\,
F_2) \in \Fln \times \Grkn \times \Fln$.  The left side of the
equation says that $V$ intersects $F_1$ $\mu$-much and also intersects
$F_2$ $\theta$-much. The right side is a union of cycles, in each of
which $V$ intersects $F_1$ $\nu$-much, has {\em no} condition directly
relating $V$ to $F_2$, and instead $F_1$ and $F_2$ are related by a
condition Poincar\'e dual to ${e''}_{\theta\mu}^\nu$. Again, the equation
is shifting the double burden on $V$ to a single burden on $V$ and a
single burden on the pair $(F_1,F_2)$.

In analogy with double Schubert calculus we feel it is appropriate to
dub the computation of the ${e''}_{\theta \mu}^\nu$ {\dfn triple
Schubert calculus}.  Note that these $e$'s are {\em not} symmetric in
$\theta$ and $\mu$, and only become so when restricted to the flag
manifold sitting diagonally in the first and third factor.

It is worth noting that triple Schubert calculus has many extensions
-- e.g. $K$-theory, replacing the Grassmannian by a flag manifold,
or using groups other than $\GLn$.\footnote{%
  It seems difficult to formulate the notion of positivity for groups
  other than $\GLn$; in \cite{MS} the roots $y_i - y_j$ of $\GLn$ are
  implicitly seen as a specialization of $y_i - z_j$, and it is
  unclear how to extend this to arbitrary root systems.  For many of
  the other extensions, no satisfactory explicit combinatorial formula
  for the structure constants is known. } 
In this way, one can view \cite{MS} as establishing a positivity
result for triple Schubert calculus on Grassmannians.%
\footnote{Much as \cite{Bu} did for K-theoretic Schubert calculus
on Grassmannians. 
In both cases, they gave a formula directly rather than an abstract reason
for positivity such as the one given in \cite{G} for equivariant
Schubert calculus.}

\subsection{An alternate interpretation: equivariant double Schubert calculus.}

For completeness, we use the connection between double Schubert calculus and
equivariant Schubert calculus discussed in subsection
\ref{double-eqvt-sec} to recast triple Schubert calculus as
``equivariant double Schubert calculus''.

Inside $\Fln \times \Grkn$, we have two interesting families of
subvarieties parameterized by $n \choose k$: the diagonal-$\GLn$ orbit
closures, and the varieties $\Fln \times X_\lambda$ corresponding to
the Schubert cycles $X_\lambda$. Both families are invariant under the
diagonal action of the torus, and so define families of equivariant
cohomology classes $\{D_\lambda\}, \{1 \tensor \tS_\mu\}$ in
$H^*_T(\Fln \times \Grkn)$.  Either family gives a $H^*_T(\Fln)$-basis
of $H^*_T(\Fln \times \Grkn)$. Therefore we can expand the product
$$ D_\theta (1\tensor \tS_\mu) = \sum_\nu {e'}_{\theta\mu}^\nu
(1\tensor \tS_\nu) $$ where the coefficients ${e'}_{\theta\mu}^\nu$
live in $H^*_T(\Fln)$.  Following a similar analysis as in subsection
\ref{double-eqvt-sec} one can show that these coefficients
${e'}_{\theta \mu}^\nu \in H^*_T(\Fln)$ refine the coefficients
${e''}_{\theta \mu}^\nu$ in triple Schubert calculus.

The results of \cite{MS} have the rather surprising consequence that
the constants ${e''}_{\theta \mu}^\nu$ and ${e'}_{\theta \mu}^\nu$ can be
lifted beyond their respective rings $H^*((\Fln)^2)$ and $H^*_T(\Fln)$
to actual double polynomials 
$e_{\theta \mu}^\nu \in H^*_{T\times T}(pt)$; this would suggest
(speaking loosely) that one should be able to replace the
cohomology ring $H^*_T(\Fln \times \Grkn)$ with the ``doubly
equivariant'' ring $H^*_{T \times T}(\Grkn)$ (using an ineffective
action of $T\times T$ on $\Grkn$).  While this can indeed be done, the
homological perspective is lost, because the classes being multiplied
are no longer representable by subvarieties.

\subsection{MS-puzzles solve the Molev-Sagan problem.}

\begin{figure}[htbp]
  \begin{center}
    \epsfig{file=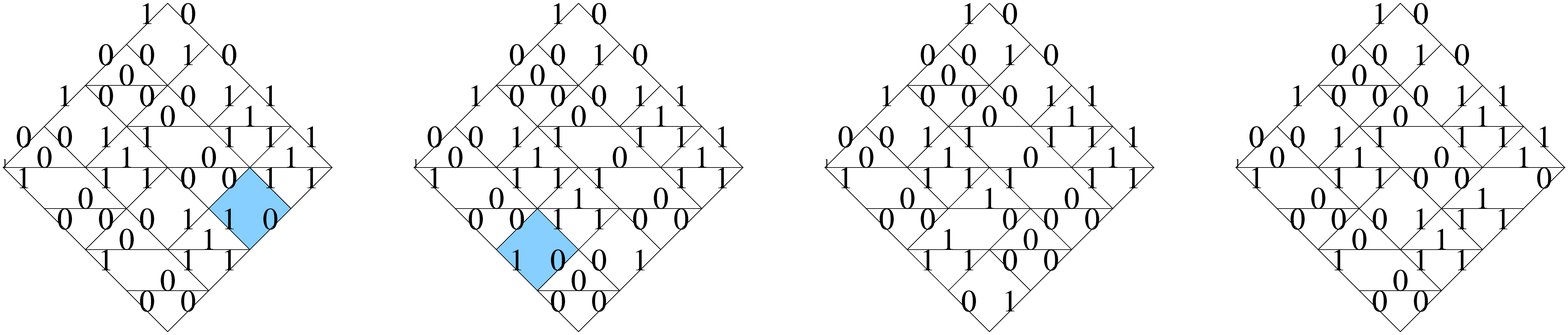,height=2in,width=6in}
    \caption{The MS-puzzles computing $ s_{0101}(x|z)\, s_{0101}(x|y) 
= ((y_3 - z_1) + (y_1 - z_2)) s_{0101}(x|y) + s_{1001}(x|y) +
s_{0110}(x|y)$.}
    \label{fig:MSexample}
  \end{center}
\end{figure}

Define an {\dfn MS-puzzle} as an arrangement of 
puzzle pieces forming a {\em diamond} of size $n$, looking something like 
a very large equivariant piece. They therefore have the usual NW and NE
boundaries, and now SW and SE boundaries as well. We will 
require the labels on the NE boundary (read clockwise) to be the
inversion-free string $\id = 0^{n-k} 1^k$.  Some examples are in
figure \ref{fig:MSexample}.

Define the {\dfn MS-weight of an equivariant piece} as $y_j - z_i$,
where $i$ measures the distance from the SE side and $j$ from the SW side%
\footnote{ This definition of weight does {\em not} reduce to the
definition we needed for (non-MS) puzzles, in the case that the
equivariant piece is in the top half.  The definition that would do
that would be $y_i-y_j$ for $i$ the distance from the SW side and $j$
the distance from the NW side.  } (starting from 1), and the {\dfn
weight of an MS-puzzle} as the product of the MS-weights of its
equivariant pieces.

\begin{Theorem*}
  The Molev-Sagan structure constant $e_{\theta\mu}^\nu$ is equal to
  the sum of the weights of the MS-puzzles with $\theta$ on the SW
  side, $\mu$ on the NW side, and $\nu$ on the SE side, all read
  bottom-to-top.  (The NE side has all $0$s, then all $1$s, read
  clockwise.)
\end{Theorem*}

We will not prove the theorem here, except to say that the Molev-Sagan
structure constants satisfy a recurrence similar to that in corollary
\ref{cor:crecurrence}, and the MS-puzzle formula can be shown to
satisfy this recurrence through a slight variant of the arguments
in section \ref{sec:cscab}.

Of course, similar recurrences were proven in \cite{MS} using barred
tableaux.  It is possible, though quite unpleasant, to establish a
weight-preserving bijection between MS-puzzles and Molev-Sagan barred
tableaux, but we will not present one here.

An example of MS-puzzles in action is in figure \ref{fig:MSexample},
demonstrating the equality
$$
\begin{array}{cccc}
 s_{0101}(x|z)\, s_{0101}(x|y) &=&& (x_1 - z_1 + x_2 - z_2) (x_1 - y_1
   + x_2 - y_2) \hfil \\ &=&& (x_2-y_1)(x_1-y_1) \hfil \\ & &+&
   ((x_1-y_2)(x_1-y_3) + (x_2-y_1)(x_1-y_3) + (x_2-y_1)(x_2-y_2) )
   \hfil \\ & &+ &((y_3 - z_1) + (y_1 - z_2)) (x_1 - y_1 + x_2 - y_2)
   \hfil \\ & =& &s_{1001}(x|y) + s_{0110}(x|y) + ((y_3 - z_1) + (y_1
   - z_2)) s_{0101}(x|y).\hfil
\end{array}
$$

\begin{figure}[htbp]
  \begin{center}
    \epsfig{file=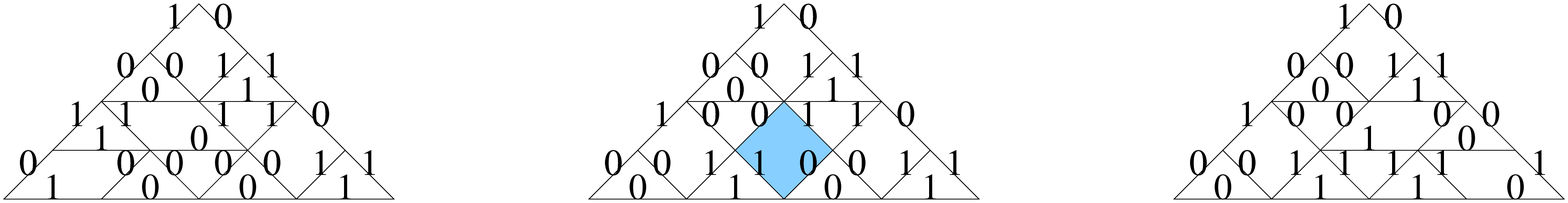,height=1in,width=5in}
    \caption{The puzzles computing $ \tS_{0101}\, \tS_{0101} 
= (y_3 - y_2) \tS_{0101} + \tS_{1001} + \tS_{0110}$.}
    \label{fig:non-MSexample}
  \end{center}
\end{figure}

It is interesting to compare this calculation to that of $\tS_{0101}^2
= (y_3-y_2)\tS_{0101} + \tS_{1001} + \tS_{0110}$ using (non-MS)
puzzles, as done in figure \ref{fig:non-MSexample}, which only uses
three puzzles.

The fact that equivariant Schubert calculus and Molev-Sagan structure
constants both reduce to ordinary Schubert calculus, in the case 
$l(\nu) = l(\lambda) + l(\mu)$, is reflected in the fact that the
two ordinary puzzles in the second calculation occur as the lower halves
of the corresponding MS-puzzles (rotated $60^\circ$).

\section{Appendix: existence of Schubert classes, 
  and the equivariant Pieri rule}

In this appendix we extend the standard {\em combinatorial} proof of
existence of Schubert classes (via divided difference operators)
to equivariant Schubert classes.\footnote{This essentially follows
Demazure's work \cite{D}, which was implicitly a calculation in
equivariant {\em K-theory} localized at the fixed points of the 
flag manifold.}
Recall that, before, 
we established the existence of equivariant Schubert classes 
by direct {\em topological} means, but this did not give a formula
for restrictions to fixed points.
As a corollary of the formula we get a direct proof 
of the equivariant Pieri rule 
(which then implies ordinary Pieri rule as a corollary).

The permutation group $S_n$ acts on $n\choose k$ and on
$H^*_T(pt) = \integers[y_1,\ldots,y_n]$ in obvious ways.
If $\alpha$ is a class, and $w\in S_n$, put these actions together
to define $w\cdot \alpha$ by
$$ (w \cdot \alpha)|_\mu := w\cdot (\alpha|_{w^{-1} \mu})$$
which is easily seen to again be a class (i.e. satisfies the GKM divisibility
conditions). 
We will care most about the case $w = s_i := (i \leftrightarrow i+1)$.

We now define the {\dfn divided difference operators} $\{\partial_i\}$. 
If $\alpha$ is a class, define $\partial_i \alpha$ by
$$ \partial_i \alpha := (\alpha - s_i\cdot \alpha) / (y_{i+1} - y_i).$$
A priori, this is just a list of rational functions. But in fact
these $\{\partial_i\}$ turn out to define endomorphisms of $H^*_T(\Grkn)$
(as a vector space):

\begin{Lemma*}
  If $\alpha$ is a class, then $\partial_i \alpha$ is also a class.
\end{Lemma*}

\begin{proof}
  From the GKM conditions we see that 
  $\partial_i \alpha \in \Oplus_{n \choose k} H^*_T(pt)$.  
  We want to know that $\partial_i \alpha$ itself satisfies the
  GKM conditions, 
  i.e. that $(\partial_i \alpha)|_\mu - (\partial_i \alpha)|_{\mu'}$
  is a multiple of $y_j - y_k$ if $\mu,\mu'$ differ in only the $j,k$
  positions. Plainly this is true if $j=i$, $k=i+1$ (or vice versa)
  since then the difference is zero.  Otherwise, the division by
  $y_{i+1} - y_i$ is irrelevant since its GCD with $y_j - y_k$ is one,
  and then the divisibility follows from the fact that $\alpha$ 
  and $s_i\cdot \alpha$ are both classes.
\end{proof}

Recall that in section \ref{ssec:Gr} we gave a topological proof of
the existence of Schubert classes (which we already knew by lemma 
\ref{lem:uniq} to be unique). The first conclusion in the following lemma 
gives a combinatorial proof, using divided difference operators,
and the second conclusion will be used in the proof of equivariant Pieri.

\begin{Lemma}\label{lem:ddopkilled}
  Fix $i \in \{1,\ldots,n-1\}$ and $\lambda\in {n\choose k}$.

  If $\lambda_i > \lambda_{i+1}$ (i.e. $s_i\cdot \lambda < \lambda$),
  then $\partial_i \tS_\lambda = \tS_{s_i \cdot\lambda}$.  
  
  If however
  $\lambda_i \leq \lambda_{i+1}$ (i.e. $s_i\cdot \lambda \geq \lambda$),
  then $\partial_i \tS_\lambda = 0$. 
\end{Lemma}

In particular, one can construct the Schubert class $\tS_\lambda$ 
by starting with the class $\tS_{w_0\cdot \id}$ (which is trivial to compute) 
and applying successive divided difference operators.

\begin{proof}
The class $\tS_\lambda$ is supported above $\lambda$, which implies
that $\partial_i \tS_\lambda$ is supported inside $\{ s_i\cdot \lambda \}
\cup \{ \mu \in {n \choose k}: l(\mu) \geq l(\lambda)\}$.  On the
other hand, from degree considerations $\partial_i \tS_\lambda$ is a
linear combination of Schubert classes of degree at most
$l(\lambda)-1$.  From these two facts and Proposition \ref{prop:basis}
we see that $\partial_i \tS_\lambda$ must vanish if $s_i\cdot \lambda \geq
\lambda$, and is an integer multiple of $\tS_{s_i\cdot \lambda}$ if $s_i\cdot
\lambda < \lambda$.  In the latter case, we can show this multiple is
1 by the straightforward computation
$$ \partial_i \tS_\lambda|_{s_i\cdot \lambda} 
  = - s_i\cdot \cdot \tS_\lambda|_\lambda / (y_{i+1} - y_i) 
  = \prod_{(j,k) \in \inv(s_i\cdot \lambda)} (y_k - y_j) 
  = \tS_{s_i\cdot \lambda}|_{s_i\cdot \lambda}.
$$
\end{proof}
%

We can now prove the equivariant Pieri rule directly.

\begin{Proposition*}[The equivariant Pieri rule]
  $$ \tS_\dv \tS_\lambda = (\tS_\dv|_\lambda) \tS_\lambda +
\sum_{\lambda': \lambda' \covers \lambda} \tS_{\lambda'}. $$
\end{Proposition*}

\begin{proof}
From Lemma \ref{structure-lem} we have
  $$ \tS_\dv \tS_\lambda = (\tS_\dv|_\lambda) \tS_\lambda +
\sum_{\lambda': \lambda' \covers \lambda} c_{\dv,\,\lambda}^{\lambda'}
\tS_{\lambda'} $$ for some integers $c_{\dv,\lambda}^{\lambda'}$; our
task is to show that $c_{\dv,\, \lambda}^{\lambda'} = 1$.

If $\lambda'\covers \lambda$, then they must differ in only two spots
  $i,i+1$, where $\lambda$ has $01$ and $\lambda'$ has $10$. Applying
  $\partial_i$, we get
$$ \partial_i(\tS_\dv \tS_\lambda) = (\tS_\dv|_\lambda) \partial_i
  \tS_\lambda + c_{\dv,\,\lambda}^{\lambda'} \partial_i \tS_{\lambda'} +
  \sum_{\mu: \mu \covers \lambda, \mu \neq \lambda'}
  c_{\dv,\,\lambda}^\mu \partial_i \tS_{\mu} $$ 
By lemma \ref{lem:ddopkilled} we have $\partial_i \tS_\lambda = 0$, hence
$s_i\cdot \tS_\lambda = \tS_\lambda$, and
$$ \partial_i(\tS_\dv \tS_\lambda) 
   = \frac{\tS_\dv \tS_\lambda - s_i\cdot(\tS_\dv \tS_\lambda)}{y_{i+1}-y_i}
   = \frac{\tS_\dv  - s_i\cdot\tS_\dv}{y_{i+1}-y_i}\,\tS_\lambda
   = (\partial_i \tS_\dv) \tS_\lambda = \tS_\lambda. $$ 
Also we have $\partial_i \tS_{\lambda'} = \tS_\lambda$ 
and $\partial_i \tS_\mu = 0$ in the above summation.  The claim
$c_{\dv, \,\lambda}^{\lambda'} = 1$ follows.
\end{proof}

Applying the forgetful map to ordinary cohomology we recover the
ordinary Pieri rule $S_\dv S_\lambda = \sum_{\lambda': \lambda'
\covers \lambda} S_{\lambda'}$.

\bibliographystyle{alpha}

\end{document}